\documentclass[a4paper, 12pt]{article}

\usepackage{latexsym,amssymb,amsmath,epsfig,amsfonts,graphicx,mathrsfs,a4wide}
\usepackage{bm}
\usepackage{pstricks,pst-node,pst-tree}
\usepackage{tikz}
\usepackage{relsize}
\usepackage{mathabx}
\usepackage[english]{babel}
\usetikzlibrary{arrows}

\usepackage{booktabs}

\newcommand{\vc}[1]{\boldsymbol{#1}}

\usepackage{dsfont}
\newcommand{\qed}{\hfill $\square$}

\newcommand{\registered}
   {{\scriptsize \ooalign{\hfil\raise0.07ex\hbox{\scriptsize \sc r}\hfil%
              \crcr\mathhexbox20D}}}

\newtheorem{theorem}{Theorem}

\newtheorem{lemma}{Lemma}

\newtheorem{corollary}{Corollary}

\newtheorem{remark}{Remark}
\newtheorem{proposition}{Proposition}
\newtheorem{conjecture}{Conjecture}

\newcommand{\nc}{\newcommand}
\nc{\ds}{\displaystyle}
\nc{\mbZ}{\mathbb Z}
\nc{\mbQ}{\mathbb Q}
\nc{\mbR}{\mathbb R}
\nc{\mbC}{\mathbb C}
\nc{\mbN}{\mathbb N}
\nc{\mbE}{\mathbb E}
\nc{\mbP}{\mathbb P}

\nc{\PH}{\emph{PH} }
\nc{\ME}{\emph{ME} }
\nc{\LST}{\emph{LST} }
\nc{\rank}{\mbox{rank\hspace{1pt}}}

\numberwithin{equation}{section}

\numberwithin{figure}{section}

\begin{document}
\title{The probabilities of extinction in a branching random walk on a strip}
\author{Peter Braunsteins\footnote{The University of Queensland, Australia. Email: p.braunsteins@uq.edu.au} \;and Sophie Hautphenne\footnote{The University of Melbourne, Australia. Email: sophiemh@unimelb.edu.au}}
\date{}
\maketitle

\begin{quote}{\small
\noindent{\bf Abstract }--
We consider a class of multitype Galton-Watson branching processes with a countably infinite type set $\mathcal{X}_d$ whose mean proge\-ny matrices have a block lower Hessenberg form. 
For these processes, the probability $\vc q(A)$  of extinction in subsets of types $A\subseteq \mathcal{X}_d$ may differ from the global extinction probability $\vc q$ and the partial extinction probability $\tilde{\vc q}$.
After deriving partial and global extinction criteria, we develop conditions for $\vc q<\vc q(A)<\tilde{\vc q}$. We then present an iterative method to compute the vector $\vc q(A)$ for any set $A$. Finally, we investigate the location of the vectors $\vc q(A)$ in the set of fixed points of the progeny generating vector.
\\

}
\end{quote}

\section{Introduction}


Multitype Galton-Watson branching processes (MGWBPs) describe the evolution of a population of independent individuals who live for a single generation and, at death, give birth to a random number of offspring that may be of various types. 
One of the primary quantities of interest in a MGWBP is the probability that the population eventually becomes empty or \emph{extinct}. 
Let $\bm{Z}_n = (Z_{n,\ell})_{\ell \in \mathcal{X}}$ record the number of type-$\ell$ individuals alive in generation $n \geq 0$ of a population whose members take types that belong to the countable set $\mathcal{X}$.
We let 
\[
q_i = \mbP[ \lim_{n \to \infty} \sum_{\ell \in \mathcal{X}} Z_{n,\ell} = 0 | \varphi_0 =i ]
\]
be the probability of \emph{global extinction} given that the population begins with a single individual of type $\varphi_0 =i$, and we refer to $\bm{q}:=(q_i)_{i \in \mathcal{X}}$ as the \emph{global extinction probability vector}.  It is well known that $\bm{q}$ is the minimal nonnegative solution of the fixed point equation $\bm{s}=\bm{G}(\bm{s})$, where $\bm{G}(\bm{s}) := (G_i(\bm{s}))_{i \in \mathcal{X}}$ records the probability generating function associated with the reproduction law of each type.

When the set $\mathcal{X}$ is finite, many of the fundamental questions concerning $\bm{q}$ are resolved in classical texts such as \cite{Har63}. 
In particular, it is well known that: \emph{(i)} $\bm{q}=\bm{1}$ if and only if the Perron-Frobenius eigenvalue of the \emph{mean progeny matrix} $M:=((\partial G_i(\vc s)/\partial s_j)|_{\vc s=\vc 1})_{i,j\in\mathcal{X}}$ is less than or equal to one, \emph{(ii)} $\bm{q}$ can be be numerically computed by repeatedly applying $\bm{G}(\cdot)$ to a vector initially comprised of zeros, and \emph{(iii)}, when $M$ is irreducible,
the set of fixed point solutions
\[
S = \{ \bm{s} \in [0,1]^{\mathcal{X}} : \bm{s} = \bm{G}(\bm{s}) \}
\]
contains at most two elements, $\bm{q}$ and $\bm{1}$. 


If we allow $\mathcal{X}$ to contain countably infinitely many types then this complicates matters considerably.
Indeed, even the definition of extinction is no longer unambiguous.
We let 
\[
\tilde{q}_i = \mbP [ \lim_{n \to \infty} Z_{n,\ell} = 0, \,  \forall \ell \in \mathcal{X} | \varphi_0 = i ]
\]
be the probability of \emph{partial extinction} given the population begins with a single individual of type $i$, and we refer to $\bm{\tilde q}:=(\tilde q_i)_{i \in \mathcal{X}}$ as the \emph{partial extinction probability vector}. 
Like $\bm{q}$, the vector $\bm{\tilde q}$ is an element of $S$.
While global extinction implies partial extinction, there may be a positive chance that every type eventually disappears from the population while the total population size grows without bound; it is then possible that $\bm{q}< \bm{\tilde q}$.
To generalise \emph{(i)}--\emph{(iii)} to the infinite type setting it is generally accepted that we should give the corresponding results for both $\bm{q}$ and $\bm{\tilde q}$. 
That is, we aim to
\emph{(i)}  derive a partial and a global extinction criterion,
\emph{(ii)} develop iterative methods to compute $\bm{q}$ and $\bm{\tilde q}$ when an algebraic expression cannot be found, and
\emph{(iii)} locate $\bm{q}$ and $\bm{\tilde q}$ in $S$.
While open questions remain, a number of authors have made progress on \emph{(i)} \cite{Braun2017,Gan10,moyal62,Spa89,Zuc11}, \emph{(ii)} \cite{Braun2018,haut12,Sag13}, \emph{(iii)} \cite{Zuc17,Braun2017,moyal62} (to name a few).

While the literature focuses on global and partial extinction, it is natural to define extinction more generally.
For $A \subseteq \mathcal{X}$ we let 
\[
q_i(A) = \mbP [ \lim_{n \to \infty} Z_{n, \ell} = 0 , \, \forall \ell \in A | \varphi_0 = i ],
\]
be the probability that the types in $A$ eventually die out given the population begins with a single individual of type $i$, and we let $\bm{q}(A):=(q_i(A))_{i \in \mathcal{X}}$ be the corresponding extinction probability vector. 
The vectors $\bm{q}(A)$ are also elements of $S$ (see Equation \eqref{FPqA}).
Such a general definition of extinction leads to redundancies. 
Indeed, in an irreducible branching process, if $A$ is finite, then $\bm{q}(A)=\bm{\tilde q}$ and $\bm{q}(\mathcal{X} \backslash A) = \bm{q}$ (see Theorem \ref{redundancy}).
However, when $A$ is infinite it is possible that $\bm{ q}<\bm{q}(A)<\bm{\tilde q}$ (see Examples 1 and 2). 
The vectors $\bm{q}(A)$ are therefore interesting in their own right. 
Apart from the recent work in \cite{Zuc14,Zuc17} which did not directly address the possibility that $\bm{q}<\bm{q}(A)<\bm{\tilde q}$, it appears that the vectors $\vc q(A)$ have received little attention in the literature.
In this more general context, Assertions \emph{(i)}--\emph{(iii)} lead to a number of natural questions: \emph{(i)} can we use $M$ to determine whether $\bm{ q}<\bm{q}(A)<\bm{\tilde q}$? \emph{(ii)} How do we compute $\bm{q}(A)$? \emph{(iii)} Can we locate the extinction probability vectors $\bm{q}(A)$ in $S$?
These questions are the primary focus of this paper.

Properties of the vectors $\bm{q}(A)$ are difficult to derive for general MGWBPs with infinitely many types. 
We therefore restrict our attention to a subclass of branching processes that is more amenable to analysis.
One possible subclass is the \emph{lower Hessenberg branching processes} considered in \cite{Braun2017}. 
In these processes, which have the typeset $\mathcal{X}=\{0,1,2,\dots\}$, the primary restriction is that \emph{type-$i$ individuals can produce offspring of type no larger than $i+1$}. 
For LHBPs, the authors of \cite{Braun2017} derive partial and global extinction criteria, and identify $\vc q$ and $\tilde{\vc q}$ respectively as the minimum and maximum of a continuum of elements in $S$. 
While tractable, these processes are too restrictive for our purposes. 
This is because in an irreducible LHBP, if $A$ is finite then $\bm{q}(A)=\bm{\tilde q}$, whereas if $A$ is infinite then $\bm{q}(A)=\bm{q}$. An irreducible LHBP therefore has at most two distinct extinction probability vectors: $\bm{q}$ and $\bm{\tilde q}$.
Here we extend  the class of LHBPs so that there may exist $A$ such that $\bm{q}< \bm{q}(A)< \bm{\tilde q}$.
We refer to processes in this extended class as \emph{block LHBPs}.
In a block LHBP, which has the typeset $\mathcal{X}_d := \{ 0,1,2, \dots \} \times \{ 1, 2, \dots, d \}$, the primary restriction is that \emph{type-$\langle i,k \rangle$ individuals can produce offspring of type $\langle j, \ell \rangle$, where $j$ is no larger than $i +1$}.
%
Following the terminology of \cite{Bol00} where random walks in a random environment on a strip (without branching) are studied, we can equivalently refer to block LHBPs as \textit{branching random walks on a strip}.

We derive a number of results for block LHBPs. We start by developing partial and global extinction criteria (Section~\ref{Part}). We then turn our attention to the more general extinction probability vectors $\bm{q}(A)$ (Section~\ref{Sets}). In particular,
\begin{itemize}
\item[\emph{(i)}] we provide sufficient conditions for $\vc q=\vc q(A)$,  $\vc q<\vc q(A)<\tilde{\vc q}$ and $\vc q(A)=\tilde{\vc q}$ (Section~\ref{Exti}),
\item[\emph{(ii)}] we develop an iterative method to compute $\vc q(A)$ for any set $A$ (Section~\ref{Comp}), and
\item[\emph{(iii)}] we make progress towards locating the vectors $\vc q(A)$ in the set $S$ (Section~\ref{secFP}).
\end{itemize}
Perhaps the most interesting part of the paper is Section~\ref{SDExample}. In this section we apply the results developed in Section~\ref{Sets} to treat an example where, by varying a single parameter, we can transition smoothly between situations where there exists one, two and four distinct extinction probability vectors. 
This example leads us to conjecture a rule for identifying which elements of $S$ correspond to an extinction probability vector $\bm{q}(A)$: we postulate that the vectors $\bm{q}(A)$ correspond to points of non-differentiability on the boundary of finite-dimensional projections of $S$ (Conjecture \ref{ConjectPH}). We conjecture that this rule extends to \textit{any} irreducible multitype Galton-Watson branching process with countably many types.


In this paper, we let $\vc 1$ and $\vc 0$ denote the infinite column vectors of $1$'s and $0$'s, respectively, and we let $\bm{1}_x$ represent the $x\times 1$ vector of $1$'s.  For any vectors $\vc x$ and $\vc y$, we write $\vc x\leq \vc y$ if $x_i\leq y_i$ for all $i$, and $\vc x< \vc y$ if $\vc x\leq  \vc y$ with $x_i<y_i$ for at least one entry $i$. Finally, we denote by $\vc e_i$ the infinite vector with all entries equal to zero, except entry $i$ which is equal to 1.

\section{Preliminaries and notation}

Consider a multitype Galton-Watson process with countably infinite type set $\mathcal{X}_d= \{  \langle k,i \rangle : k \geq 0, \, 1 \leq i \leq d \}$ for some {$ 1 \leq d < \infty$}. It will be implicitly assumed that  the types in any subset $A\subseteq\mathcal{X}_d$ are ordered lexicographically. 
 We assume that the process initially contains a single individual whose type is denoted by $\varphi_0$. 
 It then evolves according to the following rules:
\begin{itemize}
\item[\it{(i)}] each individual lives for a single generation, and
\item[ \it{(ii)}] at death gives birth to $\bm{r}=(r_{  \ell})_{{  \ell} \in \mathcal{X}_d}$ offspring, that is, $r_{\langle 0,1 \rangle }$ individuals of type $\langle 0,1 \rangle$, $r_{\langle 0,2 \rangle}$ individuals of type $\langle 0,2 \rangle$, etc., where the  vector $\bm{r}$ is chosen independently of that of all other individuals according to a probability distribution, $p_j( \cdot)$, specific to the parental type $j \in \mathcal{X}_d$.
\end{itemize}
We say that a type-$\langle k,i \rangle$ individual is in \emph{level} $k$ and \emph{phase} $i$. We partition $\mathcal{X}_d$ in two ways: by level, $\mathcal{X}_d={\bigcup}_{k \geq 0} \ell_k$, where $\ell_k= \{ \langle k,1 \rangle ,\langle k,2 \rangle , \dots, \langle k,d \rangle \}$; and by phase, $\mathcal{X}_d=\bigcup_{i=1}^d A_i$, where $A_i = \{ \langle 0,i \rangle , \langle 1,i \rangle , \dots \}$.
The primary assumption we make is that an individual in level $k$ cannot have any level $j > k+1$ offspring. In other words, the offspring vector from a level-$k$ individual belongs to the set 
\[
R_{k,d} = \left\{ \bm{r} \in \mbN_0^{\mathcal{X}_d}: \, r_{j}=0 \, \, \, \forall  j \in \bigcup^\infty_{i=k+2} \ell_i \right\}.
\]
While this assumption is made throughout, many of our results hold without it.
We refer to the resulting process as a \emph{block lower Hessenberg branching process}, or block LHBP for short. 

The branching process is defined on the Ulam-Harris space \cite[Ch. VI]{Har63}, labelled $( \Omega, \mathcal{F}, \mbP)$, as follows. Let $\mathcal{J} = \bigcup_{n \geq 0} \mathcal{J}_n$ where $\mathcal{J}_n$ describes the virtual $n$-th generation. That is, $\mathcal{J}_0= \mathcal{X}_d$, where $\varphi_0 \in \mathcal{J}_0$ specifies the type of the root, and for $n\geq 1$, $\mathcal{J}_n={  \mathcal{X}_d \times} ( \mbN \times \mathcal{X}_d \times \mbN)^n$, where $(\varphi_0;i_1, j_1, y_1; \dots; i_n, j_n, y_n)$ denotes the $i_n$-th child of type $j_n$ born to $(\varphi_0; i_1, j_1, y_1; \dots; i_{n-1}, j_{n-1}, y_{n-1})$, and $y_n$ denotes the individual's unique identification number.
\footnote{To define the branching process on the Ulam-Harris space we generally do not require that each pair $i_n, j_n$ be represented by the triple $i_n, j_n, y_n$. However, this additional index will be used in the sequel.}
Each virtual individual $I \in \mathcal{J}$ is assigned a random offspring vector $\bm{N}(I) =({N}_{\ell}(I))_{{  \ell} \in \mathcal{X}_d}$  that takes values in $R_{k,d}$ when $I$'s type belongs to level $k$, and has distribution
 ${p}_j(\cdot)$ when $I$ is of type $j$, independently of  all other individuals. The random set of individuals who appear in the population, $X = \bigcup_{n \geq 0} X_n$, is then defined recursively from the values of $\bm{N}(I)$ as follows
\[
X_0 = \{ \varphi_0 \}, \quad X_n = \{ x = ( \tilde{x}; i_n, j_n, n) \in \mathcal{J}_n : \tilde{x} \in X_{n-1}, \, i_n \leq N_{j_n}(\tilde{x}) \}.
\]
The population in generation $n$ is described by the vector $\bm{Z}_n$ with entries 
\[
Z_{n,j} = \sum_{I \in \mathcal{J}_n }  \mathds{1} ( I \in X_n , \, \, j_n = j ), \quad j \in \mathcal{X}_d.
\]
We will often refer to branching processes by their sequence of population vectors $\{ \bm{Z}_n \}_{n \geq 0}$.

From the set of probability distributions $\{ p_j (\cdot) \}_{j\in\mathcal{X}_d}$ we define the \emph{progeny generating function} $\bm{G}(\cdot) : [0,1]^{\mathcal{X}_d} \to [0,1]^{\mathcal{X}_d}$,  which contains entries
\begin{equation}
G_j ( \bm{s} ) = \sum_{ \bm{r} } p_{j}(\bm{r}) \bm{s}^{\bm{r}} = \sum_{ \bm{r} }  p_{j}(\bm{r}) \prod_{  k\in\mathcal{X}_d } s_k^{r_k},\quad j\in\mathcal{X}_d.
\end{equation}
We denote the $n$-fold composition of $\bm{G}(\cdot)$ by $\bm{G}^{(n)}(\cdot)$.
For any $n \geq 0$ and any set of types $A \subseteq \mathcal{X}_d$, let
\[
{\mathcal{E}}_n(A) = \left\{ \omega \in \Omega: \sum^\infty_{\ell=n} \sum_{i \in A} Z_{\ell,i} =0 \right\}
\]
denote the event that no individual of type $i \in A $ appear in the population from generation $n$, and let $\mathcal{E}(A)= \lim_{n \to \infty} \mathcal{E}_n(A)$ represent {the event of eventual extinction in $A$}.
{For $n\geq 0$,} we define the vector $\bm{q}_n{(A)}$ whose $i$-th element  is given by 
\[
 q_{n,i}(A) = \mbP_i ( \mathcal{E}_n(A) ),
\]
where $\mbP_i(\cdot):= \mbP(\cdot | \varphi_0=i)$. {The vector $\bm{q}_0{(A)}$, which represents the probability that no individual with type in $A$ will ever be produced, plays an important role in the sequel.} 
The vectors $\{ \bm{q}_n(A)\}_{n \geq 0}$ form a monotone non-decreasing sequence that satisfies the equation 
\[
\bm{q}_{n+1}(A) = \bm{G}(\bm{q}_n(A)),\quad n\geq 0.
\]
Consequently, by the monotone convergence theorem, each \emph{extinction probability vector} $\bm{q}(A):= \mbP(\mathcal{E}(A))$ is obtained as the limit of the sequence $\{ \bm{q}_n(A)\}$ as $n \to \infty$.
In addition, by continuity of $\bm{G}( \cdot)$, we have
\begin{equation}\label{FPqA}
\bm{q}(A)=\bm{G}( \bm{q}(A)),
\end{equation}
which implies that $\bm{q}(A)$ is an element of the set
\[
S = \{ \bm{s} \in [0,1]^{\mathcal{X}_d} : \bm{s} = \bm{G}(\bm{s}) \}.
\]
Let $T_k = \bigcup^k_{i=0} \ell_i$ be set of  types whose level is at most $k$. Following \cite{haut12} we refer to $\bigcap_{k=1}^\infty \mathcal{E}(T_k )$ and $\mathcal{E}(\mathcal{X}_d)$ as the \emph{partial} and \emph{global} extinction events, respectively, and denote their corresponding extinction probability vectors by $\bm{\tilde{q}}$ and $\bm{{q}}$. 

The mean progeny matrix $M$ is an infinite matrix whose entries are given by
\[
M(i,j) =  \left. \frac{\partial G_i ( \bm{s} )}{\partial s_j }  \right|_{\bm{s}= \bm{1}}, \quad \text{ for } i,j \in \mathcal{X}_d, 
\]
where $M(i,j)$ can be interpreted as the expected number of type-$j$ children born to a parent of type $i$. By assumption $M$ has a block lower Hessenberg structure,
\[
M=
\begin{bmatrix}
M_{00} & M_{01} & 0 & 0 & 0&  \dots \\
M_{10} & M_{11} & M_{12} & 0 & 0 & \\
M_{20} & M_{21} & M_{22} & M_{23} & 0 & \\
M_{30} & M_{31} & M_{32} & M_{33} & M_{34} & \\
\vdots & & & & & \ddots
\end{bmatrix},
\]
where for $k,l\geq 0$, $M_{kl}:=(M(i,j))_{i\in\ell_k, j\in\ell_l}$ are square matrices of order $d$.
To $M$, we associate a weighted directed graph, referred to as the \emph{mean progeny representation graph}. 
This graph has vertex set $\mathcal{X}_d$ and contains an edge from $i$ to $j$ of weight $M(i,j)$ if and only if $M(i,j)>0$. 
The branching process is said to be \emph{irreducible} if there is a path between any two  vertices in the mean progeny representation graph on $\mathcal{X}_d$.
It is well known (see for instance \cite[Proposition 4.1]{haut12}) that 
\[
\nu(M) \leq 1 \quad \Leftrightarrow \quad \bm{\tilde{q}}=\bm{1},
\]
where {$\nu(M):=\sup_{i,j}\{\limsup_n(M^n)_{ij}^{1/n}\}$} denotes the convergence norm {of $M$}.

For any $A \subseteq \mathcal{X}_d$ we define a  branching process labelled $\{ \bm{ \tilde Z}_n^{(A)} \}_{n\geq 0}$. This process, constructed on $( \Omega, \mathcal{F}, \mbP)$, is such that for any $\omega \in \Omega$,
\[
\bm{\tilde N}^{(A)}(\omega,I) = 
\begin{cases} 
\bm{N}(\omega,I), \quad & I \in A, \\
0, \quad & I \notin A,
\end{cases} 
\]{where the notation $I \in A$ means that the type of individual $I$ is in $A$.}
In other words, an outcome of $\{ \bm{ \tilde Z}_n^{(A)} \}$ corresponds to one of $\{ \bm{Z}_n \}$ in which the {individuals in $\bar{A}:=\mathcal{X}_d \backslash A$ are sterile, that is, they produce no offspring}. The process $\{ \bm{ \tilde Z}_n^{(A)} \}$  performs two roles that have parallels in the study of Markov chains on $\mathcal{X}_d$. First, $\{ ({ \tilde Z}_{n,x}^{(A)})_{x \in A} \}_{n\geq 0}$ is the branching process formed by immediately killing offspring with type in $\bar{A}$, that is, the process with the \emph{taboo} subset $\bar{A}$. Second, $(\tilde Z_{n,x}^{(T_k)})_{x \in \ell_{k+1}}$ is the vector counting the lines of descent that first enter level $k+1$ in generation $n$, that is, the vector of \emph{$n$-step first passage paths} to $\bar{T}_k$.
We let $\tilde{M}^{(A)}= (M(i,j))_{ i,j \in A}$ be the  mean progeny sub-matrix restricted to the types in $A$, and we denote by $\bm{\tilde{q}}^{(A)}$ the \emph{global} extinction probability vector of $\{ \bm{ \tilde Z}_n^{(A)} \}$. In \cite{haut12} the authors proved that $\bm{\tilde{q}}^{(T_k) }\to \bm{\tilde{q}}$ pointwise as $k \to \infty$  (where $\bm{\tilde{q}}^{(T_k) }$ was denoted as $\bm{\tilde{q}}^{(k) }$).

 Throughout this paper we make repeated use of \cite[Theorem 3.3]{Zuc14} which, for completeness, we now state and prove.
 
 \begin{theorem}\label{FCZucLemma}
For any $A \subseteq \mathcal{X}_d$ the following statements are equivalent:
\begin{itemize}
\item[(i)] $\bm{q}(A)>\bm{q}$
\item[(ii)] there exists $i \in \mathcal{X}_d$ such that $ q_{0,i}(A) > q_i$
\item[(iii)] there exists $i \in \mathcal{X}_d$ such that there is a positive chance of global survival with $\varphi_0=i$ without visiting $A$.
\end{itemize}
\end{theorem}
\textbf{Proof:} 
We obtain \emph{(ii)}$\Rightarrow$\emph{(i)} from the inequalities $\bm{q}(A) \geq \bm{q}_0(A)$ and $\bm{q}(A)\geq \bm{q}$.
The implication \emph{(i)}$\Rightarrow$\emph{(ii)} follows from the monotonicity of  $\bm{G}^{(n)}(\cdot)$ for all $n$: if  $\vc q_0(A)\leq \vc q$, then $\vc q_n(A)= \bm{G}^{(n)}(\vc q_0(A))\leq \bm{G}^{(n)}(\vc q)=\vc q$ for all $n$, which implies $\vc q(A)\leq \vc q$. 
The relations \emph{(ii)}$\Rightarrow$\emph{(iii)} and \emph{(iii)}$\Rightarrow$\emph{(i)} are immediate. \qed

\section{Partial and global extinction criteria}\label{Part}

We begin our analysis by deriving partial and global extinction criteria  for block LHBPs. These criteria extend the results in \cite[Theorem 5.1]{Braun2017}.
They are based on the sequence of $d\times d$ matrices {$\{\mathcal{M}_k\}_{k\geq 0}$} recursively defined as
\begin{eqnarray}\label{rec_mu}
\mathcal{M}_k&=&{\sum_{n\geq 0} \left(M^{(k)}\right)^n\, M_{k,k+1}},\quad k \geq 0,
\end{eqnarray}
where {
\begin{equation}\label{Mk}M^{(0)}=M_{00}, \qquad M^{(k)}:=\sum^{k}_{i=0} M_{ki}\,\mathcal{M}_{i\rightarrow k-1},\quad k\geq 1,\end{equation} with  $\mathcal{M}_{i\rightarrow j}:=\mathcal{M}_i\,\mathcal{M}_{i+1}\cdots\mathcal{M}_{j}$ for $i\leq j$. We set $\mathcal{M}_k:= \infty$} if the {series \eqref{rec_mu} does not converge, that is, if and only if $
sp({M}^{(k)})\geq 1,
$
where $sp(\cdot)$ denotes the spectral radius. If the series converges, then 
\begin{eqnarray}\label{rec_mubis}
\mathcal{M}_k&=&{[I-M^{(k)}]^{-1}\, M_{k,k+1}},
\end{eqnarray}and we can compute $\mathcal{M}_k$ recursively.
} We refer to the matrices $\mathcal{M}_k$ as \emph{step-up} matrices because of their similarity to the \emph{step-down}  probability matrices $G^{(k)}$ in \cite{Lat99} defined for level-dependent quasi-birth-and-death processes. The term ``step-up" comes from the fact that, when $\mathcal{M}_j< \infty$ for $j = 1, \dots ,k$, the matrix $\mathcal{M}_k$ records the expected number of first passage paths to $\ell_{k+1}$ that descend from a single individual in $\ell_k$, or more specifically,
\[
\mathcal{M}_k(i,j) = \mbE_{\langle k,i \rangle} \left( \sum_{n \geq 0} { \tilde Z}_{n,\langle k+1,j\rangle}^{(T_k)} \right),
\]where $\mbE_i(\cdot):= \mbE(\cdot | \varphi_0=i).$
We show this rigorously in Lemma~\ref{LemaGlob}.
In addition to the step-up matrices $\{ \mathcal{M}_k \}$, our global extinction criterion makes use of  three regularity assumptions: 
\begin{itemize}\item[(A1)] $\inf_{i} p_i (\bm{0})>0$, \item[(A2)] $\inf_{k \geq 0, \, i,j \in \{  1, \dots, d\} } p_{  \langle k,i\rangle }( 2 \bm{e}_{  \langle k,j\rangle})>0$, \item[(A3)]  $ \sup_{k\geq 0} \|  \mathcal{A}_k\|_{\infty} < \infty$, where the $d\times d^2$ matrices $\{\mathcal{A}_k\}$ satisfy the recursion
\end{itemize}
\small{\begin{eqnarray*}\mathcal{A}_0& = &\left[ I - M_{00}  \right]^{-1} \left[V_{0,00} (\mathcal{M}_0\otimes \mathcal{M}_0)+V_{0,01} (\mathcal{M}_0\otimes I)+V_{0,10} (I\otimes \mathcal{M}_0)+V_{0,11}\right],\\
\mathcal{A}_k &=& \left[ I - \sum^{k}_{i=0} M_{ki} \,\mathcal{M}_{i\rightarrow k-1} \right]^{-1}\cdot\left\{ \sum_{i=0}^{k+1} \sum_{j=0}^{k+1} V_{k,ij} \left(\mathcal{M}_{i\rightarrow k}\otimes \mathcal{M}_{j\rightarrow k}\right)\right.\\&&\qquad\qquad\qquad\left.+\sum_{i=0}^{k-1} M_{ki} \,\mathcal{M}_{i\rightarrow j-1}\, \mathcal{A}_j \left(\mathcal{M}_{j+1\rightarrow k}\otimes \mathcal{M}_{j+1\rightarrow k}\right)\right\},\; k\geq 1,\end{eqnarray*}}\normalsize the $d\times d^2$ matrices $V_{k,ij}$ have entries 
\[
V_{k,ij}(a;b,c):=\dfrac{\partial^2 G_{\langle k,a\rangle}(\vc s)}{\partial s_{\langle i,b\rangle} \partial s_{\langle j,c\rangle}} \Big|_{\vc s=\vc 1},
\] 
and $\otimes$ denotes the Kronecker product. 

Assumptions (A1) and (A2) are simple to verify in practice, whereas assumption (A3) is more challenging  but can often be verified numerically.   The matrices $\{\mathcal{M}_k\}$ and $\{\mathcal{A}_k\}$ have a physical interpretation, as we show in Lemma \ref{LemaGlob}.

\begin{theorem}\label{globpart_crit}
If $\{\bm{Z}_n \}$ is irreducible, then  
\begin{equation}\label{partial_crit}
\mathcal{M}_k   < \infty \;\; \mbox{for all } k\geq 0 \quad\Leftrightarrow\quad \tilde{\vc q}=\vc 1,
\end{equation} 
and if $\tilde{\vc q}=\vc 1$, then under Assumptions (A1)--(A3),
\begin{equation}\label{global_crit}
\sum_{k=0}^\infty  \left(\bm{1}_d^\top\, \mathcal{M}_{0\rightarrow k} \, \bm{1}_d \right)^{-1} = \infty \quad \Leftrightarrow \quad \bm{q}=\bm{1}.
\end{equation}
\end{theorem}

\tikzset{
  treenode/.style = {align=center, inner sep=0pt, text centered,
    font=\sffamily},
  arn_n/.style = {treenode, circle, black, font=\sffamily\bfseries, draw=black,
    fill=white, text width=1.2em},
  arn_r/.style = {treenode, circle, black, draw=black, 
    text width=1.2em, very thick},
  arn_x/.style = {treenode, circle, white, font=\sffamily\bfseries, draw=black,
    fill=black, text width=1.2em}
}

\begin{figure}
\begin{center}
\hspace{1.3cm}$X(\omega)$ \hspace{3cm} $f_g(X(\omega))$ \hspace{1.5cm} $f_{p,1}(X(\omega))$ \\[+2mm]
\begin{tikzpicture}[->,>=stealth',level/.style={sibling distance = 2.8cm/#1,
  level distance = 0.9cm},scale=0.9] 
\node [arn_x,minimum size=0.5cm] {\small 0,1}
    child{ node [arn_r,minimum size=0.5cm] {\small 0,2} 
            child{ node [arn_x,minimum size=0.5cm] {\small1,1}  
            	child{ node [arn_r,minimum size=0.5cm]{\small 0,2}
			child{ node [arn_r,minimum size=0.5cm]{\small 1,1}}
			child{ node [arn_r,minimum size=0.5cm]{\small 1,2}}
		} 
		child{ node [arn_x,minimum size=0.5cm]{\small 2,1}}          								
            }                          
    }
    child{ node [arn_x,minimum size=0.5cm] {\small 1,1}
            child{ node [arn_r,minimum size=0.5cm] {\small 0,2}
            	child{ node [arn_r,minimum size=0.5cm]{\small 1,1}
			child{ node [arn_r,minimum size=0.5cm]{\small 1,2}}
			child{ node [arn_r,minimum size=0.5cm]{\small 1,2}}
		}
            } 													            
            child{ node [arn_x,minimum size=0.5cm] {\small 2,1}
							child{ node [arn_r,minimum size=0.5cm] {\small 1,1}}
							child{ node [arn_x,minimum size=0.5cm] {\small 3,1}
								child{ node [arn_x,minimum size=0.5cm]{\small 4,1}}
								child{ node [arn_x,minimum size=0.5cm]{\small 4,2}}
							}
            }
		}
;  
\end{tikzpicture}
\begin{tikzpicture}[->,>=stealth',level/.style={sibling distance = 2.8cm/#1,
  level distance = 0.9cm},scale=0.9] 
\node [arn_x,minimum size=0.5cm] {\small 0,1}
    child{ node [arn_x,minimum size=0.5cm] {\small 1,1}
    	child{ node [arn_x,minimum size=0.5cm]{\small 2,1}}
    }
    child{ node [arn_x,minimum size=0.5cm] {\small 1,1}
    	child{ node [arn_x,minimum size=0.5cm] {\small 2,1}
		child{ node [arn_x,minimum size=0.5cm]{\small 3,1}
			child{ node [arn_x,minimum size=0.5cm]{\small 4,1}}
			child{ node [arn_x,minimum size=0.5cm]{\small 4,2}}
		} 
   	 }
    }
;  
\end{tikzpicture}
\begin{tikzpicture}[->,>=stealth',level/.style={sibling distance = .8cm/#1,
  level distance = 0.9cm},scale=0.9] 
\node [arn_x,minimum size=0.5cm] {\small 1,1}
	child{ node [arn_r,minimum size=0.5cm]{\small 1,1}}
	child{ node [arn_r,minimum size=0.5cm]{\small 1,2}
		child{ node [white] {\small 1,1}
		[white] child{ node [white]{\scalebox{1}{1}}
			[white] child{ node [white]{\scalebox{1}{1}}}
		}                          
    }
	}
;  
\end{tikzpicture}
\, 
\begin{tikzpicture}[->,>=stealth',level/.style={sibling distance = 1.5cm/#1,
  level distance = 0.9cm},scale=0.9] 
\node [arn_x,minimum size=0.5cm] {\small 1,1}
	child{ node [arn_r,minimum size=0.5cm]{\small 1,1}
		child{ node [arn_r,minimum size=0.5cm]{\small 1,2}}
		child{ node [arn_r,minimum size=0.5cm]{\small 1,2}
				child{ node [white] {\small 1,1}
		[white] child{ node [white]{\scalebox{1}{1}}}                          
    }
		}
	}
;  
\end{tikzpicture}
\end{center}
\caption{\label{SDembV}A visualisation of the embedded processes for a specific $\omega \in \Omega$.   Nodes correspond to individuals, the first digit denoting their level and the second their phase.
}
\label{RF2}
\end{figure}
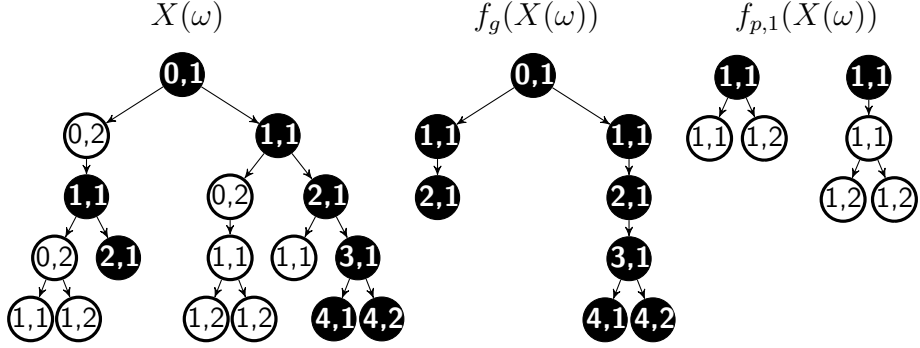

\noindent
Before turning to the proof of Theorem \ref{globpart_crit},   we introduce the branching process $\{ \bm{Y}_k \}$ embedded within $\{ \bm{Z}_n \}$, whose sample paths are constructed from those of $\{ \bm{Z}_n \}$ as follows:
we define a function $f_g(\cdot) : \mathcal{J} \to \mathcal{J}$ that takes each line of descent $(\varphi_0; i_1, j_1,  y_1;  \dots; i_n, j_n,   y_n)$ and deletes each  triple $(i_k, j_k,  y_k)$ whose level is not strictly larger than that of all its ancestors.
For each $\omega \in \Omega$ the family tree of $\{\bm{Y}_k \}$ is then given by $f_g(X(\omega))$; see the middle tree in Figure \ref{SDembV} for an example.\footnote{This is where the third entry in the triple $(i_k,j_k,y_k)$ is important. Suppose we removed it. Then for $\omega$ illustrated in Figure \ref{SDembV} we have $f_g(\langle 0,1 \rangle ;1, \langle 0,2 \rangle ;1,\langle 1,1\rangle )=f_g(\langle 0,1 \rangle ;1,\langle  1,1\rangle )=(\langle 0,1 \rangle ;1, \langle 1,1 \rangle )$, causing both individuals in generation one of $f_g(X(\omega))$ to have the same label. This makes the construction of the lineages in the next generation unclear.
}
Observe that generation $k$ of $\{ \bm{Y}_k \}$ contains individuals in level $k$ only. 
 Specifically, it contains the individuals that are the first to enter level $k$ in their line of descent. To avoid confusion we take the convention that $\{ \bm{Y}_k \}$ starts at the generation corresponding to the level of the initial type in $\{\vc Z_n\}$. 
The embedded process $\{ \bm{Y}_k \}$ evolves as a $d$-type Galton-Watson process whose offspring distributions vary deterministically with the generation: an individual's phase corresponds to its type and an individual's level corresponds to its generation. The process $\{ \bm{Y}_k\}$ is therefore a multitype Galton-Watson process in a varying environment (see for instance \cite{Jon97}). In addition, for the reasons laid out in \cite[Section 3]{Braun2017}, individuals in $\{ \bm{Y}_k \}$ may have an infinite number of offspring; {in this case, we} say that $\{ \bm{Y} _k \}$ is \emph{explosive}. According to the arguments in \cite[Corollary 1]{Braun2017},
\begin{equation}\label{MultLink1}
\mathcal{E}(\mathcal{X}_d ) \stackrel{a.s.}{=} \left\{ \lim_{k \to \infty} |\bm{Y}_k | = 0 \right\},
\end{equation}
and 
\begin{equation}\label{MultLink2}
\bigcap_{k=0}^\infty \mathcal{E}( T_k) \stackrel{a.s.}{=}  \left\{ | \bm{Y}_k | < \infty, \, \forall \, k\geq 0 \right\},
 \end{equation}
where  $| \bm{Y}_k |$ denotes the total size of generation $k$.
In other words, $\{ \bm{Z} _n \}$ experiences global extinction if and only if $\{ \bm{Y}_k \}$ experiences extinction, and $\{ \bm{Z} _n \}$ experiences local survival (avoids partial extinction) if and only if $\{ \bm{Y}_k \}$ experiences explosion. 
This enables us to evaluate whether partial or global extinction occurs in $\{ \bm{Z} _n \}$ simply by observing $\{ \bm{Y}_k \}$. 

We denote the progeny generating function of $\{ \bm{Y}_k \}$ at generation $k$ by $\bm{g}_k(\bm{{s}}^{(d)})=(g_{k,i}(\bm{{s}}^{(d)}))_{1\leq i\leq d}$, where $\bm{{s}}^{(d)} \in [0,1]^d$. For $\ell\leq k$ we let $\bm{g}_{\ell\rightarrow k}(\bm{{s}}^{(d)}):= \bm{g}_\ell \circ  \bm{g}_{\ell+1}\circ\cdots \circ\bm{g}_k(\bm{{s}}^{(d)})$.

\begin{lemma}\label{implicit}
For any $k\geq 0$, the progeny generating function $\bm{g}_k(\cdot)$ satisfies
\begin{equation}\label{impliciteq}
{g}_{k,i}(\bm{{s}}^{(d)})=G_{ \langle k,i \rangle}(\bm{g}_{0 \to k}(\bm{{s}}^{(d)}), \bm{g}_{1 \to k}(\bm{{s}}^{(d)}), \dots , \bm{g}_k(\bm{{s}}^{(d)}),\bm{{s}}^{(d)}).
\end{equation} 
\end{lemma}
\textbf{Proof:} The proof follows the same conditioning argument as that of Lemma 3 in \cite{Braun2017} but in the multitype setting. \qed

 We now show that the matrices $\{\mathcal{M}_k\}$ and $\{\mathcal{A}_k\}$ correspond to the first and second factorial moment progeny matrices in $\{\vc Y_k\}$.
 
 \begin{lemma}\label{LemaGlob}
If $\bm{ \tilde q}=\bm{1}$, then for any $k\geq 0$ ,
\[
\mathcal{M}_k(i,j) = \left. \frac{\partial {g}_{k,i}(\bm{s}^{(d)})}{\partial s^{(d)}_{j}} \right|_{\bm{s}^{(d)}=\bm{1}},  \quad \mbox{ and } \quad \mathcal{A}_k(i;j,l) =  \left. \frac{\partial^2 {g}_{k,i}(\bm{s}^{(d)})}{\partial s^{(d)}_{j} \partial s^{(d)}_{l}} \right|_{\bm{s}^{(d)}=\bm{1}}.
\]

\end{lemma} 
\textbf{Proof:} By \eqref{MultLink2} and the assumption $\bm{\tilde{q}}=\bm{1}$, we have $|\bm{Y}_k|<\infty$ almost surely for all $k$. Thus, $\bm{g}_{\ell\to k} (\bm{1} )=\bm{1}$ for all $\ell\leq k$ and $k\geq 0$. The statement then follows by successive differentiations of \eqref{impliciteq}. \qed
 
 \medskip

\noindent
\textbf{Proof of Theorem \ref{globpart_crit}:} 
By Lemma \ref{LemaGlob}, assertion \eqref{global_crit} follows from \eqref{MultLink1} and \cite[Theorem 2.3]{Dolgo2017}.

To obtain \eqref{partial_crit},
{we embed a second process} in $\{ \bm{Z}_n \}$, this time with the mean progeny matrix $M^{(k)}$ defined in \eqref{Mk}. To do this we introduce a function $f_{p,k}(\cdot): \mathcal{J} \to \mathcal{J}$ that takes a (possibly infinite) line of descent $(\varphi_0;i_1,j_1,{  y_1}; i_2 , j_2,  y_2,\dots)$, and operates in two stages: first, it deletes the descendants of all triples $(i_\ell, j_{\ell},   y_\ell)$ whose level is strictly larger than $k$, to obtain the corresponding line of descent in $\{ \bm{\tilde{Z}}_n^{(T_{k})} \}$; and second, it deletes all remaining triples whose level differs from $k$ to obtain the \emph{restriction} (see \cite[p118]{Lat99}) of $\{ \bm{\tilde{Z}}_n^{(T_{k})} \}$ to level $k$.  
When the function $f_{p,k}(\cdot)$ is applied to a random tree $X$, the result is a random tree which evolves as a $d$-type Galton Watson process; see the right tree in Figure \ref{SDembV} as an example. 
In addition, if  $\mathcal{M}_j<\infty$ for all $j=0,\dots, k-1$, the mean progeny matrix of this embedded process is indeed given by $M^{(k)}$. By irreducibility,  with probability 1 this embedded process endures extinction if and only if $\{ \bm{\tilde Z}_n^{(T_{k})}\}$ does as well. Invoking the extinction criterion for finite-type processes (see for instance \cite[Ch. II, Theorem 7.1]{Har63}), $\bm{\tilde q }^{(T_{k})} =\bm{1}$ if and only if  $sp(M^{(k)}) \leq 1$.
Thus, if $sp(M^{(k)})<1$ for all $k\geq 0$, then $\bm{\tilde q }^{(T_{k})} =\bm{1}$ for all $k$, and according to \cite[Theorem A.1]{Braun2017} we then have $\bm{ \tilde q}=\bm{1}$. Similarly, if there exists $k$ such that $sp(M^{(k)}) >1$ then $\bm{\tilde q }^{(T_{k})} <\bm{1}$ and $\bm{ \tilde q}\leq \bm{\tilde{q}}^{(T_{k})}<\bm{1}$. 
Finally, if there exists $k$ such that $sp(M^{(k)})=1$, then by irreducibility there exists a path from level $k$ to itself via a maximum level $\ell >k$ in the mean progeny representation graph of $M$, which again leads to $\bm{ \tilde q}\leq \bm{\tilde{q}}^{(T_{\ell})}<\bm{1}$. \qed

\section{Extinction in sets of types}\label{Sets}

\subsection{Extinction criteria}\label{Exti}

We now shift our attention to the more general extinction probability vectors $\bm{q}(A)$, in particular, we investigate how to determine when $\bm{q}(A)$ differs from $\bm{q}$ and $\bm{ \tilde q}$. 
We begin with a general result that allows us to use $q_{0,i}(A)$, the probability that a type-$i$ individual has no descendants in $A$, to compare extinction probability vectors.
\begin{theorem}\label{redundancy}
Let $A,B \subseteq \mathcal{X}_d$. If $\sup_{i \in B} {q}_{0,i}(A)<1$ then \mbox{$\bm{q}(A) \leq \bm{q}(B)$}.
\end{theorem} 
\textbf{Proof:} 
Let $\mathcal{F}_n$ denote the history of the process up to generation $n$. 
By L\'evy's 0-1 law, for any  fixed $\ell \geq 0$, {\begin{equation}\label{conv_prob}\mbP( \mathcal{E}_\ell(A) | \mathcal{F}_n ) \to { \mathds{1}}(\mathcal{E}_\ell(A))\quad \mbox{as $n \to \infty$}\end{equation}}on a subset $\Omega^*_{\ell}$ of the sample space that has probability 1. Let $\Omega^* = \bigcap_{\ell \geq 0} \Omega^*_{\ell}$. For any outcome $\omega \in \bar{\mathcal{E}}(B) \cap \Omega^*$ (such that $\{ \bm{Z}_n(\omega) \}$ contains  individuals with types in $B$ for infinitely many $n$), {we have} $\mbP( \mathcal{E}_\ell(A) | \mathcal{F}_n)(\omega) < 1 - \varepsilon$ for infinitely many $n$, and for some $\varepsilon>0$. Thus, by {\eqref{conv_prob}}, $\mathds{1}(\mathcal{E}_\ell(A))(\omega)<1-\varepsilon$, that is, $\omega \in \bar{\mathcal{E}}_{\ell}(A)$. Since this holds for all $\ell$, $\omega \in \bigcup_{\ell \geq 0} \bar{\mathcal{E}}_\ell(A) = \bar{\mathcal{E}}(A)$. Hence  $\mathcal{E}(A)\cap \bar{\mathcal{E}}(B) \subseteq \bar{\Omega}^*$, leading to 
\[
\mbP_i(\mathcal{E}(A)) = \mbP_i\left(\mathcal{E}(A) \cap\mathcal{E}(B) \right) + \mbP_i\left(\mathcal{E}(A) \cap \bar{\mathcal{E}}(B)\right) \leq \mbP_i( \mathcal{E}(B))    
\]
for any $i \in \mathcal{X}_d$.
\qed

\begin{corollary}\label{FinitePartial} Let $A \subseteq \mathcal{X}_d$. 
If $\{\bm{Z}_n \}$ is irreducible then $ \bm{ q}(A)\leq  \bm{\tilde{q}}$, and if  in addition $|A| < \infty$ then
 $\bm{ q}(A)= \bm{\tilde{q}}$.
\end{corollary}
\textbf{Proof:} We first show that if $|A|<\infty$, then $\bm{ q}(A)= \bm{\tilde{q}}$.
 By irreducibility, the condition of Theorem \ref{redundancy} is satisfied for any finite sets $A$ and $B$. Thus, letting $B=T_k$, we have $\bm{q}(A) = \bm{q}(T_k)$ for all $k \geq 0$. Because $\mathcal{E}(T_{k+1}) \subseteq \mathcal{E}(T_{k})$, by the monotone convergence theorem, 
\[
 \bm{q}(A)=\lim_{k \to \infty} \bm{q}(T_k)=\mbP\left(\lim_{k \to \infty}\mathcal{E}(T_{k})\right)= \bm{\tilde q}.
\]
Now, for any $A \subseteq \mathcal{X}_d$ (not necessarily finite) and $i \in A$ we have $\bm{q}(A) \leq \bm{q}(\{ i \} ),$ and by what precedes, $\bm{q}(\{ i \} )= \bm{\tilde q}$, therefore $ \bm{ q}(A)\leq  \bm{\tilde{q}}$. \qed

Given Corollary \ref{FinitePartial} we will focus on extinction in infinite sets $A$. In particular,  we shall consider sets $A$ belonging to the sigma algebra generated by the phase partition $\{ A_i \}$,
which we denote by $\sigma(A_1, \dots, A_d)$. 
As we will see, even with  just two phases ($d=2$), it is possible for a process to survive in phase one, $A_1$, while enduring extinction in phase two, $A_2$, and vice versa. A concrete example is provided in Section \ref{SDExample}. Nonetheless, the following result states that if the phases are sufficiently intertwined,  then the probability of extinction in any set $A \in\sigma(A_1, \dots, A_d)$ coincides with the global extinction probability.

\begin{corollary}\label{InfiniteGlobal}
If $\sup_{\ell \in A_i} {q}_{0,\ell}(A_j)< 1$ for all $i,j \in \{ 1, \dots, d \}$ then $\bm{q}(A)=\bm{q}$ for any $A \in \sigma(A_1, \dots, A_d)$. 
\end{corollary}
\textbf{Proof:} Since $d<\infty$, $\sup_{\ell \in A_i} {q}_{0,\ell}(A_j)< 1$ for all $i,j \in \{ 1, \dots, d \}$ implies $\sup_{i \in \mathcal{X}_d} {q}_{0,i}(A)< 1$ for any $A \in \sigma(A_1, \dots, A_d)$. The statement then follows from Theorem \ref{redundancy}.
\qed

 Corollaries \ref{FinitePartial} and \ref{InfiniteGlobal} indicate that, under quite general conditions, $\bm{q}(A)=\bm{\tilde{q}}$ if $|A|<\infty$ and $\bm{q}(A)=\bm{{q}}$ if $|A|=\infty$,  the same as in the single-phase LHBP analysed in \cite{Braun2017}. So, when do we have  $\bm{q} < \bm{q}(A)< \bm{\tilde q}$? We begin with a necessary condition, which follows from Theorem \ref{FCZucLemma}.
 
 \begin{corollary}\label{NecDif}
If $\bm{q} < \bm{q}(A)$ then 
\begin{equation}\label{NecC}
\bm{\tilde{q}}^{( \bar{A})} <\bm{1}.
\end{equation}
\end{corollary}

Corollary \ref{NecDif} states that to have $\bm{q} < \bm{q}(A)$, it must be possible for $\{ \bm{Z}_n \}$ to survive in the types $\bar A$ without any outside assistance from the types in $A$. 
 To verify \eqref{NecC}, we observe that when $A \in \sigma(A_1, \dots,A_d)$, $\tilde M^{(\bar A)}$ is block lower Hessenberg; we can then compute the sequence $\{ \tilde{\mathcal{M}}_k ^{( \bar A)} \}_{ k \geq 0} $ using \eqref{rec_mu} with $\tilde{M}^{(\bar A)}$ substituted for $M$, and apply Theorem \ref{globpart_crit}. The matrices $\{ \tilde{\mathcal{M}}_k^{(\bar A)} \}$ are also a fundamental ingredient in Theorem \ref{OCrit}. In preparation for this theorem, for each level $k\geq 0$, we let $\bar{A}(k)=\bar A \cap \ell_k$, and we define
\begin{itemize}\item
 the column
vector $\bm{t}_k^{({\bar{A}})}=(t_{k,i}^{(\bar{A})})_{i\in \bar{A}(k)}$, where
\[
t_{k,i}^{(\bar{A})} := \sum_{ j \in A} M({ i} , j)
\] 
is the
expected total number of  direct descendants in $A$ from a parent of type $\langle k,i \rangle \in \bar{A}$, and
\item the matrix $\tilde F_k^{(\bar{A})}=(\tilde F_k^{(\bar{A})}(i,j))_{i, j \in \bar{A}(k)}$, where  $\tilde F_k^{(\bar{A})}(i,i):=1, $ and where for $i\neq j$,
\[
\tilde F_k^{(\bar{A})}(i,j) := M^{(\bar{A})}(i , j)+\displaystyle{\sum_{\substack{
k\geq 1\\ i_1,i_2, \dots, i_{k} \neq j}}}\tilde M^{(\bar{A})}(i , i_1)\tilde M^{(\bar{A})}(i_1,i_2) \cdots \tilde M^{(\bar{A})}(i_k,j)
\]
is the weighted sum of first passage paths from $i$ to $j$ in level $k$ in the mean progeny representation graph of $\tilde M^{(A)}$. 
\end{itemize}
We also let $\tilde{\mathcal{M}}^{(\bar{A})}_{0\rightarrow k-1}:=\mathcal{M}_0^{(\bar A)}\mathcal{M}_1^{(\bar A)}  \dots \mathcal{M}_{k-1}^{(\bar A)}$, and $v$ be the number of phases in $\bar A$ so that $v = |\bar A(k)|$ for all $k$.

\begin{theorem}\label{OCrit}
Let $A \in \sigma(A_1, \dots, A_d)$, and assume $\bm{\tilde{q}}^{(\bar{A})}<\bm{1}$ and $\nu(\tilde{M}^{(\bar{A})}) <1$.  If, in addition,
\begin{itemize}
\item[(A)] $\sum^{\infty}_{k=0} (\bm{1}_v^\top \bm{t}_k^{(\bar A)}) \,\tilde{ \mathcal{M}}^{(\bar{A})}_{0\rightarrow k-1} \bm{1}_v < \infty, $
and 
\item[(B)] there exists $K<\infty$ such that $\tilde F_k^{(\bar A)} \leq   K\, \vc 1_v\cdot\vc1_v^\top$ for all $k \geq 0$,
\end{itemize}
 then $\bm{q} < \bm{q}(A)$ and $\bm{q}(\bar A) <\bm{\tilde q}$.
\end{theorem}

\noindent\textbf{Proof:} 
 We first demonstrate that, under the conditions of the theorem, the expected number of sterile individuals produced over the lifetime of $\{ \bm{\tilde{Z}}_n^{(\bar{A})} \}$ (those with type in $A$) is finite.  Without loss of generality we assume that the process starts with an individual of type $i\in \bar A(0)$.
  Let $\tilde M^{(\bar{A},n)}( i ,j)$ denote the $(i,j)$th entry of the $n$th power of $\tilde M^{(\bar{A})}$.
The expected number of sterile types produced throughout the lifetime of $\{ \bm{\tilde{Z}}_n^{(\bar{A})} \}$ is then given by 
\begin{align}\label{BLHBconsd}
\mbE_{i} \left( \sum_{n=0}^\infty \sum_{l \in A} \tilde{Z}_{n,l}^{({  \bar{A}})} \right) & = \sum_{k=0}^\infty \sum_{j\in \bar{A}(k)} \sum_{n=0}^\infty \tilde M^{(\bar{A},n)}( i ,j) \,t_{k,j}^{(\bar{A})}.\end{align} 
Observe that for any $k \geq 0$, $ i \in\bar{A}(0)$, and $j\in \bar{A}(k)$, we have 
\begin{align}
\sum^\infty_{n =0} \tilde M^{(\bar{A},n)}( i ,j) &=\sum_{l \in \bar{A}(k)} \tilde{\mathcal{M}}^{(\bar{A})}_{0\rightarrow k-1} (i,l)\, F_k^{(\bar{A})}(l,j)\,\sum^\infty_{n =0} \tilde M^{(\bar{A},n)}( j ,j) \nonumber\\
&\leq\left( \frac{1}{1-\nu(\tilde{M}^{(\bar{A})}) }\right) \left[\tilde{\mathcal{M}}^{(\bar{A})}_{0\rightarrow k-1}\, F_k^{(\bar{A})}\right](i,j), 
\label{SinglePe2}
\end{align}
where \eqref{SinglePe2} follows from \cite[Theorem A4]{Sen06}. 
By \eqref{BLHBconsd}, \eqref{SinglePe2}, and the assumptions of the theorem, we then have in matrix form
\begin{align*}
\mbE_{\bar{A}(0)} \left( \sum_{n=0}^\infty \sum_{l \in A} \tilde{Z}_{n,l}^{({  \bar{A}})}\right) &\leq \left( \frac{1}{1-\nu(\tilde{M}^{(\bar{A})}) }\right) \sum_{k=0}^\infty \tilde{\mathcal{M}}^{(\bar{A})}_{0\rightarrow k-1}\, \tilde F_k^{(\bar A)} \bm{t}_k^{(\bar A)} \\
&\leq  \left( \frac{K}{1-\nu(\tilde{M}^{(\bar{A})}) }\right) \sum_{k=0}^\infty ( \bm{1}_v^\top \bm{t}_k^{(\bar A)} )\,\tilde{\mathcal{M}}^{(\bar{A})}_{0\rightarrow k-1} \bm{1}_v \\
&<\infty.
\end{align*}
Because this expectation is finite, with probability 1 there exists a generation $n$ after which a sterile type in $A$ never appears in the population. Thus, under the assumption $\bm{\tilde{q}}^{(\bar{A})}<\bm{1}$,  there exists a type $i \in \bar A$ such that starting from $i$ there is a positive chance of global survival without entering the set $A$. By Theorem \ref {FCZucLemma} we then have $\bm{q} < \bm{q}(A)$. 
In addition, by the assumption $\nu(\tilde{M}^{(\bar{A})}) <1$, if $\{ \bm{Z}_n \}$ survives in $\bar A$ but not in $A$, then it becomes partially extinct with probability one, leading to $\bm{q}(\bar A) <\bm{\tilde q}$.
\qed

If for some set $B \in \sigma(A_1, \dots , A_d)$ the conditions of Theorem~\ref{OCrit} hold with both $A=B$ and $A=\bar B$, then $\bm{q}< \bm{q}(B) < \bm{\tilde q}$ and $\bm{q}< \bm{q}(\bar B) < \bm{\tilde q}$.
Condition (A) of Theorem~\ref{OCrit} can be verified easily if there are only finitely many edges between $A$ and $\bar{A}$ in the mean progeny representation graph of ${M}$, because in that case there are only finitely many values of $k$ such that $\bm{t}_k^{(\bar A)}$ is non-zero.
Condition (B) of Theorem \ref{OCrit} is of a more technical nature.   It holds for example if $\bar A$ contains a single phase, or if the phases in $\bar A$ are sufficiently intertwined, or if there is some symmetry between the phases.
The next lemma formalises this.

\begin{lemma}
If $\nu(\tilde M^{(\bar A)})<1$ then each of the following conditions are sufficient for Condition (B):
\begin{itemize}
\item[(B1)]  $\bar A=A_i$ for some $i\in\{1,\ldots,d\}$;
\item[(B2)] There exists $\varepsilon >0$ such that $\tilde F_k^{(\bar A)}( i,j) > \varepsilon$ for all $i,j \in \bar A(k)$ and $k \geq 0$;
\item[(B3)] For any $k,\ell \geq 0$, $i\in \bar A(k)$ and $j\in \bar A(\ell)$, we have $\tilde M^{(\bar A)}( i,j ) = \tilde M^{(\bar A)}( j,i )$.
\end{itemize}
\end{lemma}
\textbf{Proof:}  The sufficiency of (B1) is trivial since then $\tilde F^{(\bar A)}=1$. The sufficiency of (B2) and (B3) follows from the fact that when $\nu(\tilde M^{(\bar A)}) < 1$, for any $k\geq 0$ and $i,j\in\bar A(k)$, $\tilde F^{(\bar A)}(i,j)\tilde F^{(\bar A)}(j,i)$ is bounded above by the weighted sum of first \textit{return} paths from $i$ to $i$ in the mean progeny representation graph of $M$, which is strictly less than 1.
\qed

In the specific case where $\{ \bm{Z}_n \}$ is singular, {that is, each individual produces exactly one offspring with probability one}, the process survives with probability 1 ($\vc q=\vc 0$), and the process $\{Z_n\}_{n\geq 0}$, where ${Z}_n:=i\Leftrightarrow  {Z}_{n,i}>0$, corresponds to an irreducible Markov chain on the state space $\mathcal{X}_d$. The arguments in the proof of Theorem \ref{OCrit} then lead to the following corollary, which can be viewed as the algorithmic complement to the more theoretical result of \cite[Corollary 8]{Doo17}. 

\begin{corollary}
If $| \bm{Z}_n| = 1$ a.s. for all $n \geq 0$, then 
\begin{equation}\label{DooCor}
\bm{q}(A)>\bm{0} \quad \Leftrightarrow \quad \lim_{k \to \infty} \bm{1}_v^\top \left(  \prod_{j=k}^\infty \tilde{\mathcal{M}}_j^{(\bar{A})}\right)\bm{1}_v >0.
\end{equation}
In addition, if   $\sum_{x\in\bar A(k+1)} \tilde M^{(\bar A)}( \langle k,i \rangle, x) > 0$
 for all $\langle k , i \rangle \in \bar A$, then the right hand side of \eqref{DooCor} may be replaced by $\bm{1}_v^\top \left(  \prod_{j=0}^\infty \tilde{\mathcal{M}}_j^{(\bar{A})}\right)\bm{1}_v>0$.
\end{corollary}

\noindent\textbf{Proof.} By Theorem \ref{FCZucLemma}, $\bm{q}(A)>\bm{0} $ if and only if there exists $i\in\bar{A}$ such that the probability that the chain $\{Z_n\}$ never visits $A$ starting from $i$ is strictly positive, which is equivalent to the right hand side of \eqref{DooCor}. The additional condition  $\sum_{x\in\bar A(k+1)} \tilde M^{(\bar A)}( \langle k,i \rangle, x) > 0$ for all $\langle k , i \rangle \in \bar A$ ensures that there is no null factor in the product   $\prod_{j=0}^\infty  \tilde{\mathcal{M}}_j^{(\bar{A})}$.\qed

\subsection{Computational methods}\label{Comp}

Given the existence of extinction probability vectors $\bm{q}(A)$ different from $\bm{q}$ and $\bm{\tilde{q}}$, we now  develop a method of computing them. 

For any {$k, \ell\geq -1$}, define the finite-type branching process $\{ \bm{Z}^{(k,\ell)}_n(A)\}$ on the same probability space as $\{ \bm{Z}_n \}$, with progeny generating vector $\bm{G}^{(k,\ell)}(A)(\bm{s})$ such that
\[
{G}_i^{(k,\ell)}(A)(\bm{s})= \begin{cases}
0, \quad & \text{if } i \in A \cap \bar T_k \\
1, \quad & \text{if } i \in \bar A \cap \bar T_{\ell} \\
G_i( \bm{s}) \quad & \text{otherwise},
\end{cases}
\]
and denote by $\bm{q}^{(k,\ell)}(A)$ its extinction probability vector.
In other words, $\bm{q}^{(k,\ell)}(A)$ is the probability that the branching process $\{ \bm{ \tilde Z}_n^{(A \cup T_\ell)} \}$ becomes extinct before producing a type in $A \cap \bar T_k$ {(with $T_{-1}=\varnothing$ and $\bar T_{-1}=\mathcal{X}_d$)}. 
For any $k$ and $\ell$, the vector $\bm{q}^{(k,\ell)}(A)$ can then be computed using established techniques for finite-type branching processes.

\begin{theorem}\label{COD} If $\{ \bm{Z}_n \}$ is irreducible then
\[
\lim_{k \to \infty} \lim_{\ell \to \infty} \bm{q}^{(k,\ell)}(A)= \bm{q}(A).
\]
\end{theorem}
\textbf{Proof:} By Theorem A.1 of \cite{Braun2017}, for any fixed value of $k$,  $\lim_{\ell \to \infty} \bm{q}^{(k,\ell)}(A)$ is the partial extinction probability of the original process modified so that types in $A \cap \bar T_k$ are immortal. 
Let {\[ N(A):= \inf \left\{ N\geq 0 : \sum_{n=N}^\infty \sum_{i \in A} Z_{n,i}> 0 \right\} \]}be the last generation at which a type in $A$ appears in the population, and
 \begin{equation}\label{tauk}\tau_k(A):= \inf \left\{k\geq 0: \sum_{n=0}^k \sum_{i \in A \bigcap \bar T_k} Z_{n,i}> 0 \right\}\end{equation} be the first generation at which a type in $A \cap \bar T_k$ appears in the population.
By Corollary \ref{FinitePartial} and the fact that $|A \cap T_k |< \infty$ for all $k$, we have, for all $i\in\mathcal{X}_d$,
 \begin{align*}
 \lim_{\ell \to \infty} {q}_i^{(k,\ell)}(A)&=\mbP_i \left({ \{ N(A) < \tau_k(A) \} \bigcap \left\{ \bigcap_{u=1}^\infty \mathcal{E}\left(\bar A \cap T_u\right)\right\}}\right) \\
 &=\mbP_i \left( \{ N(A) < \tau_k(A) \} \right).
 \end{align*}By the monotone convergence theorem we then have
 \begin{eqnarray*}
\lim_{k \to \infty} \mbP_i ( \{ N(A) < \tau_k(A) \}) &=& \mbP_i ( \{ N(A) < \lim_{k \to \infty} \tau_k(A) \}) \\
&=&  \mbP_i (N(A) < \infty ) \\
&= & \mbP_i (\mathcal{E}(A)),
\end{eqnarray*}which conlcudes the proof.
\qed

Following the arguments in the proof of Theorem \ref{COD}, we obtain a method to compute the probability $\bm{q}_0(A)$ that no individuals with type in $A$ will ever be produced:
\begin{corollary} If $\{ \bm{Z}_n \}$ is irreducible then 
\[ \lim_{\ell \to \infty} \bm{q}^{(-1,\ell)}(A)= \bm{q}_0(A).\]
\end{corollary}

While Theorem \ref{COD} may be applied in a general setting, it requires both $k$ and $\ell$ to be increased to infinity separately. 
  From a computational perspective it would be more   efficient to set $\ell = k$ and let  them both increase to infinity together. We now derive a sufficient condition ensuring convergence of the resulting sequence.

\begin{theorem}\label{COD2}
 If $\sup_{i \in A } \tilde{q}_i^{(A)}<1$ then 
 \[ \lim_{k \to \infty} \bm{q}^{(k,k)}(A) = \bm{q}(A). \]
\end{theorem}
\textbf{Proof:} 
First, suppose that $\{ \bm{Z}_n \}$ becomes extinct in the set $A$. In this case, there exists $K$ such that $\tau_k(A)=\infty$ for all $k >K$, where $\tau_k(A)$ is defined in \eqref{tauk}. In addition, by Theorem \ref{redundancy} there is almost sure partial extinction. This implies that, for $k>K$ there is almost sure global extinction in $\{ \bm{Z}_n^{(k,k)}(A) \}$, leading to $\bm{q}(A) \leq \liminf_k \bm{q}^{(k,k)}(A)$.

Now suppose that $\{ \bm{Z}_n \}$ survives in the set $A$. At any generation $n$ consider the daughter processes of each individual in $( Z_{n,i})_{i \in A}$ truncated so that all types in $\bar A$ have no offspring. 
Note that if one of these truncated daughter processes survives globally, then there exists $K$ such that $Z^{(k,k)}(A)$ survives globally for all $k>K$. 
This is because for these values of $k$ an immortal type $A \cap \bar T_k$ must eventually be born into the population. 
Let $D$ be the event that, throughout the life of $\{ \bm{Z}_n \}$, there exists an individual that has a truncated daughter process which survives globally.
By assumption, there exists $\varepsilon >0$ such that, whenever $( Z_{n,i})_{i \in A}$ is non-empty, we have $\mbP (D | \mathcal{F}_n ) > \varepsilon$. Because $\{ \bm{Z}_n \}$ survives in the set $A$, $( Z_{n,i})_{i \in A}$ is non-empty for infinitely many values of $n$, therefore, following the same arguments as in the proof of Theorem \ref{redundancy}, the event $D$ occurs with probability 1. This then implies $\bm{q}(A) \geq \limsup_k  \bm{q}^{(k,k)}(A)$. \qed

\medskip
To understand why we impose the   sufficient condition  $\sup_{i \in A } \tilde{q}^{(A)}_i<1$ in Theorem \ref{COD2}, consider an example with two phases where this condition is not satisfied. Assume $\bm{G}( \bm{s})$ contains entries
\[
G_{\langle k, i \rangle}( \bm{s} ) = \begin{cases}
s_{\langle k+1, 1 \rangle} s_{\langle 0, 2 \rangle}, \quad & \langle k, i \rangle \in A_1 \\
\frac{  k+1}{  k+2}s_{\langle k+1,2 \rangle} + \frac{1}{  k+2}, \quad & \langle k, i \rangle \in A_2.
\end{cases}
\]
In this case $q_{\langle 0 ,1\rangle}(A_2)= 0$ but $q^{(k,k)}_{\langle 0 ,1 \rangle}(A_2)   = (1- \frac{1}{k+2})^{k+1} \to e^{-1}$. While this is a reducible example, it highlights a pathology that can also occur in the irreducible setting.

\subsection{Fixed Points}\label{secFP}

We now briefly turn our attention to the set $S$, which contains the extinction probability vectors $\bm{q}(A)$. To allow the results on (single phase) LHBPs provided in \cite[Section 4]{Braun2017} to be applied directly, 
we introduce the concept of a block LHBP that is \emph{locally isomorphic} to a (single phase) LHBP. 
Roughly speaking, this is a block LHBP which takes the distribution of a LHBP when we sum the number of individuals in each level (i.e $ (\sum_{j \in \ell_k} Z_{n,j})_{k \geq 0}\stackrel{d}{=} \bm{\hat Z}_n,$ where $\{ \bm{\hat Z}_n \}$ is a LHBP). 
More specifically, $\{\bm{Z}_n \}$ is locally isomorphic to a LHBP $\{ \bm{\hat Z}_n \}$ if for each level $k \geq 0$, there exists a probability distribution $\hat p_k(\cdot) : R_{k,1} \to [0,1]$ such that for any $\bm{u} \in R_{k,1}$ and $j \in \ell_k$,
\begin{equation}\label{ILHBP}
\hat p_k(\bm{u}) = \sum_{ \substack{ \bm{r}\in R_{k,d} \,s.t.\\ \sum_{x \in \ell_i} r_x=u_i \, \forall i}} p_j( \bm{r}).
\end{equation}
We define the projection of $S$ onto the $k$th level,
\[
S_{k} := \{ \bm{u} \in [0,1]^{\ell_k} : \exists \bm{s} \in S \text{ with } (s_j)_{j\in\ell_k} = \bm{u} \}.
\]
The next proposition is the irreducible block analogue of \cite[Theorem 4.1]{Braun2017}.

\begin{proposition}\label{IsoFP}
Suppose $\{ \bm{Z}_n \}$ is locally isomorphic to a LHBP that satisfies the conditions of \cite[Theorem 4.1]{Braun2017}. Then 
\[
\bm{q}=\min S,\qquad\bm{\tilde{q}} = \sup S \backslash \{ \bm{1} \},
\]
and
\begin{equation}\label{affine}
(x q_{\langle k,1\rangle} + (1 - x) \tilde{q}_{\langle k, 1 \rangle }) \bm{1}_{d} \in S_{k}, \quad \mbox{ for all } x \in [0,1],
\end{equation}where $q_{\langle k,1\rangle}=\cdots=q_{\langle k,d\rangle}$, and $\tilde{q}_{\langle k, 1 \rangle }=\cdots=\tilde{q}_{\langle k, d \rangle }.$
\end{proposition}
\textbf{Proof:} 
Let $\bm{s}= [0,1]^{\mathcal{X}_1}$ and suppose $\bm{s} = \bm{ \hat G}( \bm{s} )$. Then for any $k \geq 0$ and any $j \in \ell_k$,
\begin{align*}
{G}_j(s_0 \bm{1}_{d}, s_1 \bm{1}_{d}, \dots) &= \sum_{\bm{u} \in R_{k,d}} (s_0 \bm{1}_{d}, s_1 \bm{1}_{d}, \dots)^{\bm{u}} p_j(\bm{u}) \\
&= \sum_{\bm{u} \in R_{k,1}} \bm{s}^{\bm{u}} \left(  \sum_{\bm{v} : u_i = \sum_{x \in \ell_i} v_x, \, \forall i} p_j( \bm{v}) \right) \\
&=  \sum_{\bm{u} \in R_{k,1}} \bm{s}^{\bm{u}} \hat p_k( \bm{u} ) =s_k.
\end{align*}
Therefore,
\[
(s_0 \bm{1}_{d}, s_1 \bm{1}_{d}, \dots) =  \bm{G}(s_0 \bm{1}_{d}, s_1 \bm{1}_{d}, \dots),
\]
and the result follows as a direct consequence of \cite[Theorem 4.1]{Braun2017}. \qed

\begin{remark}
Blackwell \cite{Black1955} demonstrates that, if $\{\bm{Z}_n \}$ is singular (so $\bm{q}=\bm{0}$) and $\bm{\tilde q}=\bm{1}$, then there are \emph{two} distinct extinction probability vectors if and only if every bounded solution to $\bm{s}=\bm{G}(\bm{s})$ is of the form $\bm{s} =c\bm{1}$ for $c \geq 0$. In the context of this section, this means that if $\{\bm{Z}_n \}$ singular then $\bm{q}(A) \in \{ \bm{q}, \bm{\tilde q} \}$ for all $A \subseteq \mathcal{X}$ if and only if $S_k$ is entirely made up of the linear segment identified in Proposition \ref{IsoFP}. The next example suggests that this criterion does not generalise to non-singular branching processes.
\end{remark}

\section{Illustrations}\label{SDExample}

We now illustrate our main theorems through two examples, and motivate some open questions.

\begin{figure}
\centering
\begin{tikzpicture}[scale=0.85]

\tikzset{vertex/.style = {shape=circle,draw,minimum size=1.3em}}
\tikzset{edge/.style = {->,> = latex'}}

%

\node[vertex, minimum size=.9cm] (0) at  (0,0) {\footnotesize$0 ,1$};
\node[vertex, minimum size=.9cm] (2) at  (3.5,0) {\footnotesize $1, 1 $};
\node[vertex, minimum size=.9cm] (4) at  (7,0) {\footnotesize $2,1$};
\node[vertex, minimum size=.9cm] (6) at  (10.5,0) {\footnotesize $3,1$};

\node[vertex, minimum size=.9cm] (1) at  (0,2.5) {\footnotesize $ 0, 2 $};
\node[vertex, minimum size=.9cm] (3) at  (3.5,2.5) {\footnotesize $1,2$};
\node[vertex, minimum size=.9cm] (5) at  (7,2.5) {\footnotesize $2,2$};
\node[vertex, minimum size=.9cm] (7) at  (10.5,2.5) {\footnotesize $3,2$};

\node at (13,2.5) {\small $\dots$};
\node at (13,0) {\small $\dots$};



%

\draw[edge,above] (0) to[bend left=15] node {\footnotesize $c$ } (2);
\draw[edge,above] (2) to[bend left=15] node {\footnotesize $c $ } (4);
\draw[edge,above] (4) to[bend left=15] node {\footnotesize $c $ } (6);

\draw[edge,below] (2) to[bend left=15] node {\footnotesize $a$ } (0);
\draw[edge,below] (4) to[bend left=15] node {\footnotesize $a$ } (2);
\draw[edge,below] (6) to[bend left=15] node {\footnotesize $a$ } (4);

\draw[edge,above] (1) to[bend left=15] node {\footnotesize $c$ } (3);
\draw[edge,above] (3) to[bend left=15] node {\footnotesize $c$} (5);
\draw[edge,above] (5) to[bend left=15] node {\footnotesize $c$} (7);

\draw[edge,below] (3) to[bend left=15] node {\footnotesize $a$ } (1);
\draw[edge,below] (5) to[bend left=15] node {\footnotesize $a$} (3);
\draw[edge,below] (7) to[bend left=15] node {\footnotesize $a$} (5);


\draw[edge,left] (0) to[bend left=20] node {\footnotesize $y$ } (1);
\draw[edge,left] (2) to[bend left=20] node {\footnotesize $y/x$} (3);
\draw[edge,left] (4) to[bend left=20] node {\footnotesize $y/x^{2}$} (5);
\draw[edge,left] (6) to[bend left=20] node {\footnotesize $y/x^{3}$} (7);

\draw[edge,right] (1) to[bend left=20] node {\footnotesize $y$ } (0);
\draw[edge,right] (3) to[bend left=20] node {\footnotesize $y/x$} (2);
\draw[edge,right] (5) to[bend left=20] node {\footnotesize $y/x^{2}$} (4);
\draw[edge,right] (7) to[bend left=20] node {\footnotesize $y/x^{3}$} (6);

\draw[edge,below] (0) to [out=-110,in=-70,looseness=4.5] node {\footnotesize $b$} (0);
\draw[edge,below] (2) to [out=-110,in=-70,looseness=4.5] node {\footnotesize $b$} (2);
\draw[edge,below] (4) to [out=-110,in=-70,looseness=4.5] node {\footnotesize $b$} (4);
\draw[edge,below] (6) to [out=-110,in=-70,looseness=4.5] node {\footnotesize $b$} (6);

\draw[edge,above] (1) to [out=110,in=70,looseness=4.5] node {\footnotesize $b$} (1);
\draw[edge,above] (3) to [out=110,in=70,looseness=4.5] node {\footnotesize $b$} (3);
\draw[edge,above] (5) to [out=110,in=70,looseness=4.5] node {\footnotesize $b$} (5);
\draw[edge,above] (7) to [out=110,in=70,looseness=4.5] node {\footnotesize $b$} (7);


%

\draw[above] (7) to[bend left=7.5, pos=1] node {\footnotesize $c$} (12.25,2.825);
\draw[edge,below] (12.25,2.175) to[bend left=8,pos=0] node {\footnotesize $a$}  (7);

\draw[above] (6) to[bend left=7.5, pos=1] node {\footnotesize $c$} (12.25,0.325);
\draw[edge,below] (12.25,-0.325) to[bend left=8,pos=0] node {\footnotesize $a$}  (6);

\end{tikzpicture}
\caption{\label{Exam1MPRG}The mean progeny representation graph corresponding to Example 1.}
\end{figure}
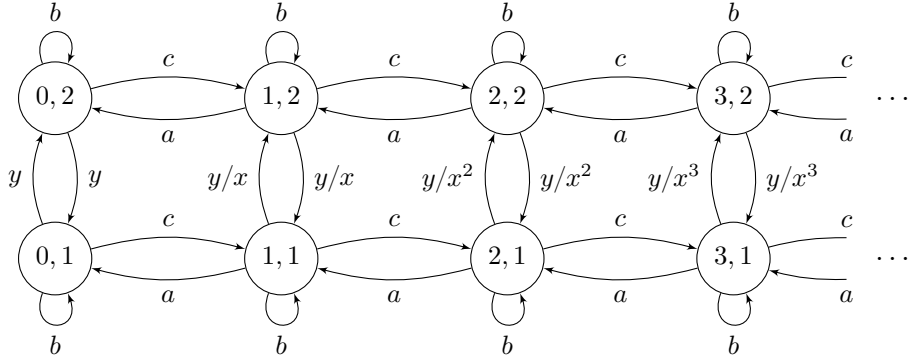

\paragraph{Example 1}
We consider a two-phase ($d=2$) block LHBP whose progeny generating vector $\bm{G}( \bm{s} )$ contains entries 
\[
G_{\langle k,i \rangle}(\bm{s})= \begin{cases}
\frac{b}{u}s_{\langle 0,i \rangle}^u+\frac{c}{u}s_{\langle 1,i \rangle}^u +\frac{y}{u}s_{\langle 0, 3-i \rangle}^u +1-\frac{b+c+y}{u}, \quad & \langle k, i \rangle \in \ell_0 \\[0.5em]
\frac{a}{u}s_{\langle k-1,i \rangle}^u+\frac{b}{u}s_{\langle k,i \rangle}^u+\frac{c}{u}s_{\langle k+1,i \rangle}^u  & \\\qquad\;\;\;\;\;\;\,+\frac{y}{u x^k}s_{\langle k, 3-i \rangle}^u +1-\frac{a+b+c+y x^{-k}}{u}, \, \,   &  \langle k, i \rangle \in \bar \ell_0 ,
\end{cases}
\]
where $a,b,c,y>0$, $x\geq 1$, and $u=\lceil a+b+c+y +1 \rceil$.
The corresponding mean progeny representation graph is illustrated in Figure \ref{Exam1MPRG}. 

In this example, 
the   processes $\{ \bm{\tilde Z}_n^{(A_1)} \}$ and $\{ \bm{\tilde Z}_n^{(A_2)} \}$ restricted to each phase form stochastically identical nearest-neighbour branching random walks with an absorbing barrier.
Individuals within the same level give birth to each other with a probability that decays geometrically at rate $x$ with the individual's level. 
  We can also verify that the process is locally isomorphic to a LHBP $\{ \bm{ \hat Z}_n \}$.   The next proposition highlights the contrasting asymptotic behaviours of the branching process as a function of the decay rate $x$.
\begin{proposition}\label{SDExam}
Suppose $b+2\sqrt{ac}<1$ and 
\[ \mu:= \left(1-b - \sqrt{(1-b)^2-4ac} \right)/2a >1. \] 
We have
\begin{itemize}
\item[(i)] if $x=1$ and $b+y+2\sqrt{ac} \leq 1$, then $\bm{q}=\bm{q}(A_1)=\bm{q}(A_2)<\bm{\tilde{q}}=\bm{1}$;
\item[(ii)] if $x=1$ and $b+y+2\sqrt{ac} > 1$, then $\bm{q}=\bm{q}(A_1)=\bm{q}(A_2)=\bm{\tilde{q}}<\bm{1}$;
\item[(iii)] if $x>1$, then $\bm{q}<\bm{\tilde{q}}$;
\item[(iv)] if $x>\mu$, then   $\bm{q}< \bm{q}(A_1)<\bm{\tilde{q}}$ and $\bm{q}<  \bm{q}(A_2)<\bm{\tilde{q}}$. 
\end{itemize} 
\end{proposition} 
\textbf{Proof:} \textit{(i) and (ii).} Suppose $x=1$.   By  Corollary \ref{InfiniteGlobal}, $\bm{q}=\bm{q}(A_1)=\bm{q}(A_2)$, and by \cite[Corollary 3]{Braun2017}, $\bm{q}<\vc1$.
Note that there is partial (global) extinction in $\{ \bm{Z}_n \}$ if any only if there is partial (global) extinction in   its local isomorphism $\{ \bm{ \hat Z}_n \}$. Denote the mean progeny matrix of $\{ \bm{ \hat Z}_n \}$ by $\hat M$. This is a tridiagonal matrix with entries $\hat M(i,i-1) =a \bm{1}\{ i \geq 1\}$, $\hat M(i,i)=b+y$, $\hat M(i,i+1)=c$ for $i \geq 0$, and $0$ otherwise. By \cite[Proposition 5.1]{haut12}, $\nu(\hat M)=b+y+2\sqrt{ac}$, which means
\[
\bm{ \tilde q}  = \bm{1} \quad \Leftrightarrow \quad b+y+2\sqrt{ac} \leq 1
\]
yielding \emph{(i)}. Observe that, when $x=1$, $\{ \bm{ \hat Z}_n \}$ is a branching random walk with an absorbing barrier, which implies that $\tilde q_{\langle k , i \rangle }$ is decreasing in $k$. Consequently, when $b+ y + 2\sqrt{ac} >1$ the entries of $\bm{ \tilde q}$ are uniformly bounded away from 1. By \cite[Lemma 3.3]{moyal62}, $S$ contains only one such element, which, when combined with the fact that $\bm{q}\leq \bm{\tilde q}$, yields \emph{(ii)}.


\textit{(iii).} Suppose $x>1$.
Let $\wideparen M^{(\bar T_k)} = (\wideparen M^{(\bar T_k)}(i,j) )_{i,j \geq 1}$, where $\wideparen M^{(\bar T_k)}(i,j) := \hat M (\langle i + k,1 \rangle, \langle j+k,1 \rangle)$,     is the mean progeny matrix of the process $\{ \bm{ \hat Z}_n \}$ taboo on $T_k$   (which we denote by $\{ \bm{ \wideparen Z}^{(\bar T_k)}_n \}$), with entries relabelled for convenience.  It is such that $\wideparen M^{(T_k)}(i,i)=b+y/x^{k+i}$, $\wideparen M^{(T_k)}(i,i+1)=c$ for $i\geq 1$, and $\wideparen M^{(T_k)}(i,i-1)=a$ for $i\geq 2$. We have $b\leq b+y/x^{k+i}\leq b+y/x^{k+1}$ for all $k,i$; therefore, by definition of the convergence norm, and by \cite[Proposition 5.1]{haut12}, $b + 2\sqrt{ac}\leq \nu( \wideparen M^{(\bar T_k)} )\leq b+y/x^{k+1} + 2\sqrt{ac}$ for all $k$, leading to
%
\[
\lim_{k \to \infty} \nu( \wideparen M^{(\bar T_k)} ) = b + 2\sqrt{ac}<1.
\]
In addition, by \cite[Corollary 3]{Braun2017}, for any $k \geq 0$ and initial type $i \in \bar T_k$, $\{ \bm{ \wideparen Z}^{(\bar T_k)}_n \}$ has a positive chance of global survival. The assertion then follows through direct application of \cite[Theorem 7.1]{Braun2017}.

\textit{(iv). }Suppose $x>\mu$. We have $t_k^{(A_1)}=y/x^k$, and by \cite[Lemma 9]{Braun2017} $\tilde{\mathcal{M}}_k^{(A_1)} \to \mu$. Thus, 
\[
\lim_{k \to \infty} \left(\tilde{\mathcal{M}}^{(A_1)}_{0\rightarrow k-1}\,
t_k^{(A_1)}\,\right)^{1/k} =\mu/x<1.
\]
The root test for convergence then gives $  \sum_{k=0}^\infty  \tilde{\mathcal{M}}^{(A_1)}_{0\rightarrow k-1}\,
t_k^{(A_1)}<\infty$. 
By \cite[Corollary 3]{Braun2017} and \cite[Proposition 5.1]{haut12} we have $\bm{ \tilde q}^{(A_1)}< \bm{1}$ and $\nu( \tilde{M}^{(A_1)})=b + 2\sqrt{ac}<1$. 
  By Theorem \ref{OCrit} we have $\bm{q}<\bm{q}(A_2)$ and $\bm{q}(A_1)<\bm{\tilde q}$. The result then follows by repeating the same arguments with $A_2$ in place of $A_1$.
 \qed
 
 \begin{figure}
\centering
\includegraphics[width=10cm]{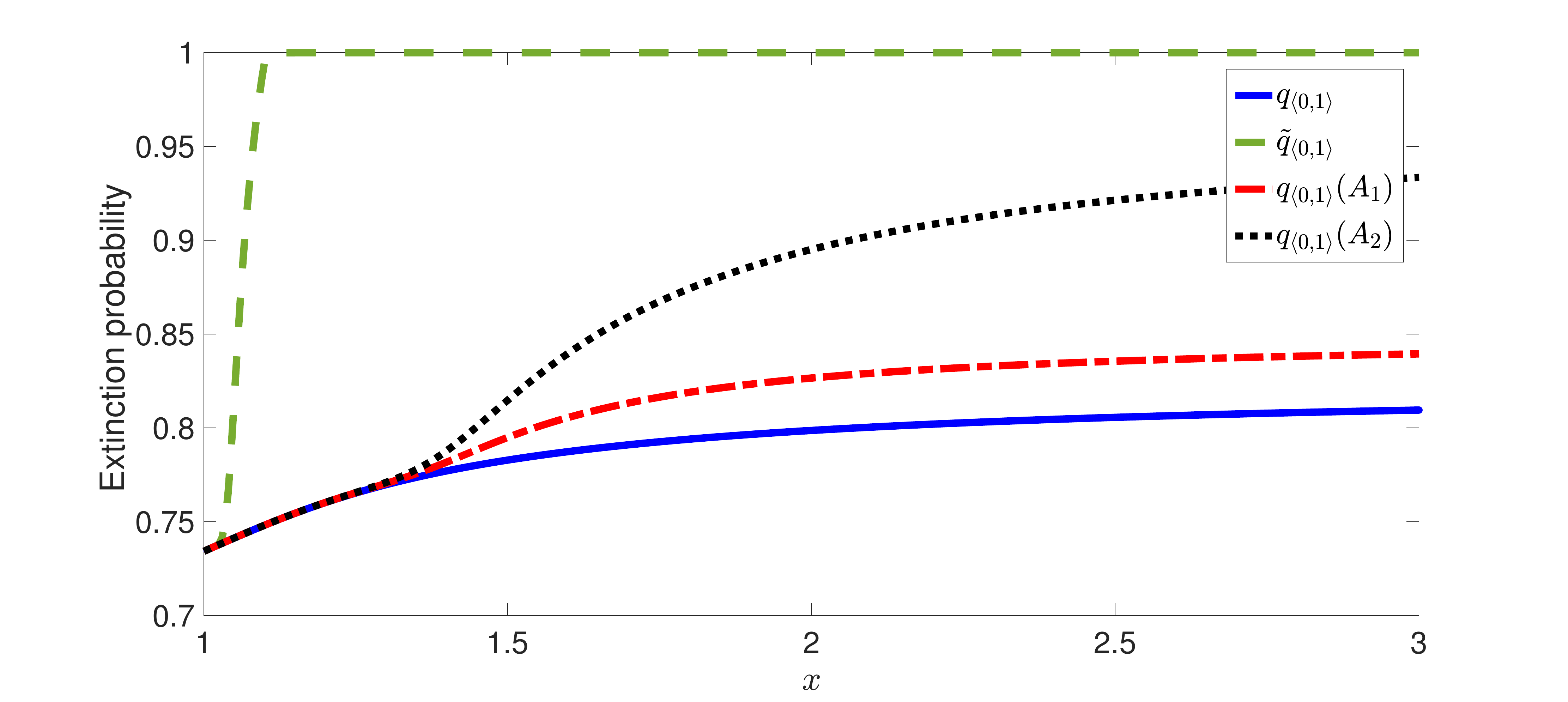}
\vspace{2mm}

\includegraphics[width=3.8cm]{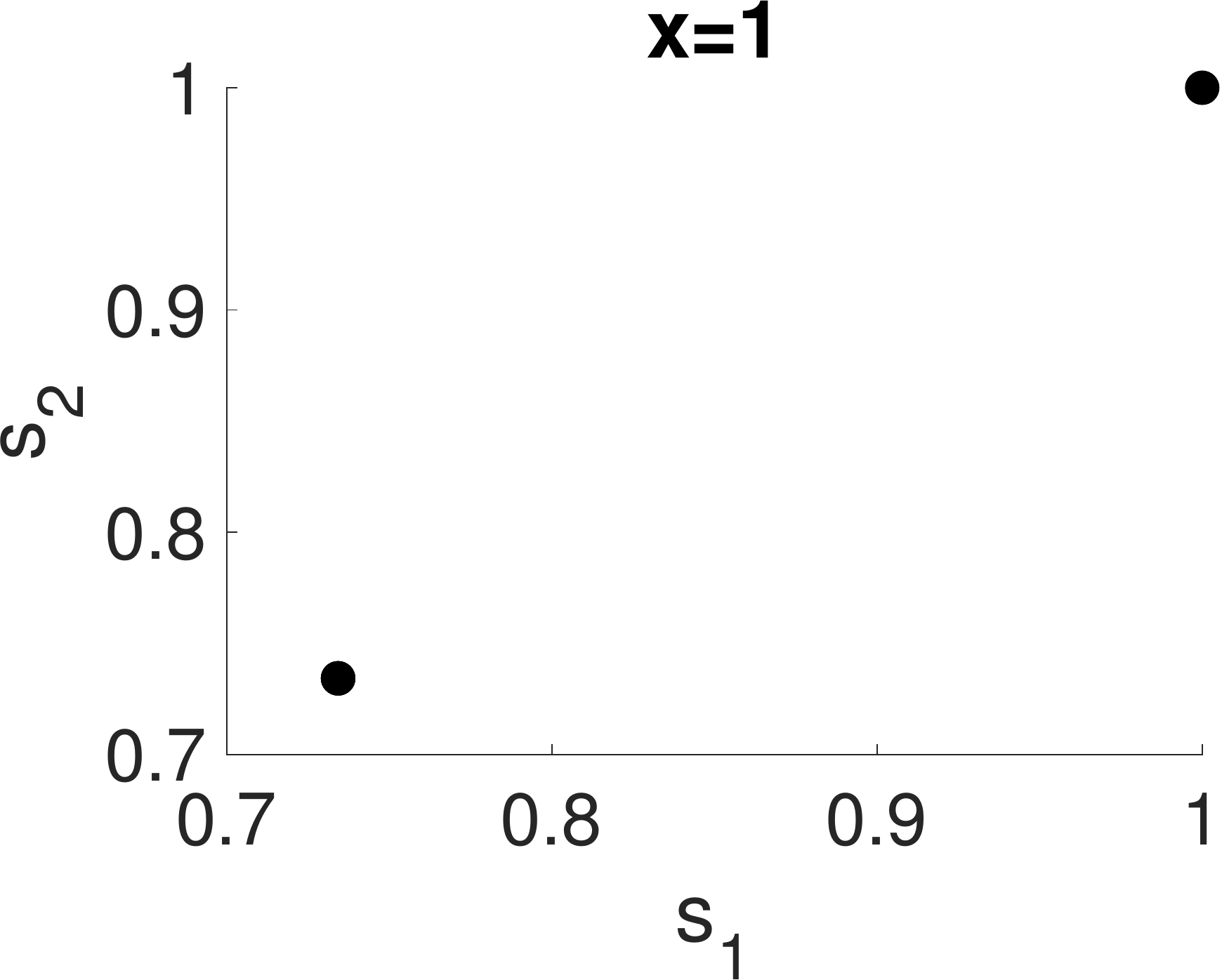}
\includegraphics[width=3.8cm]{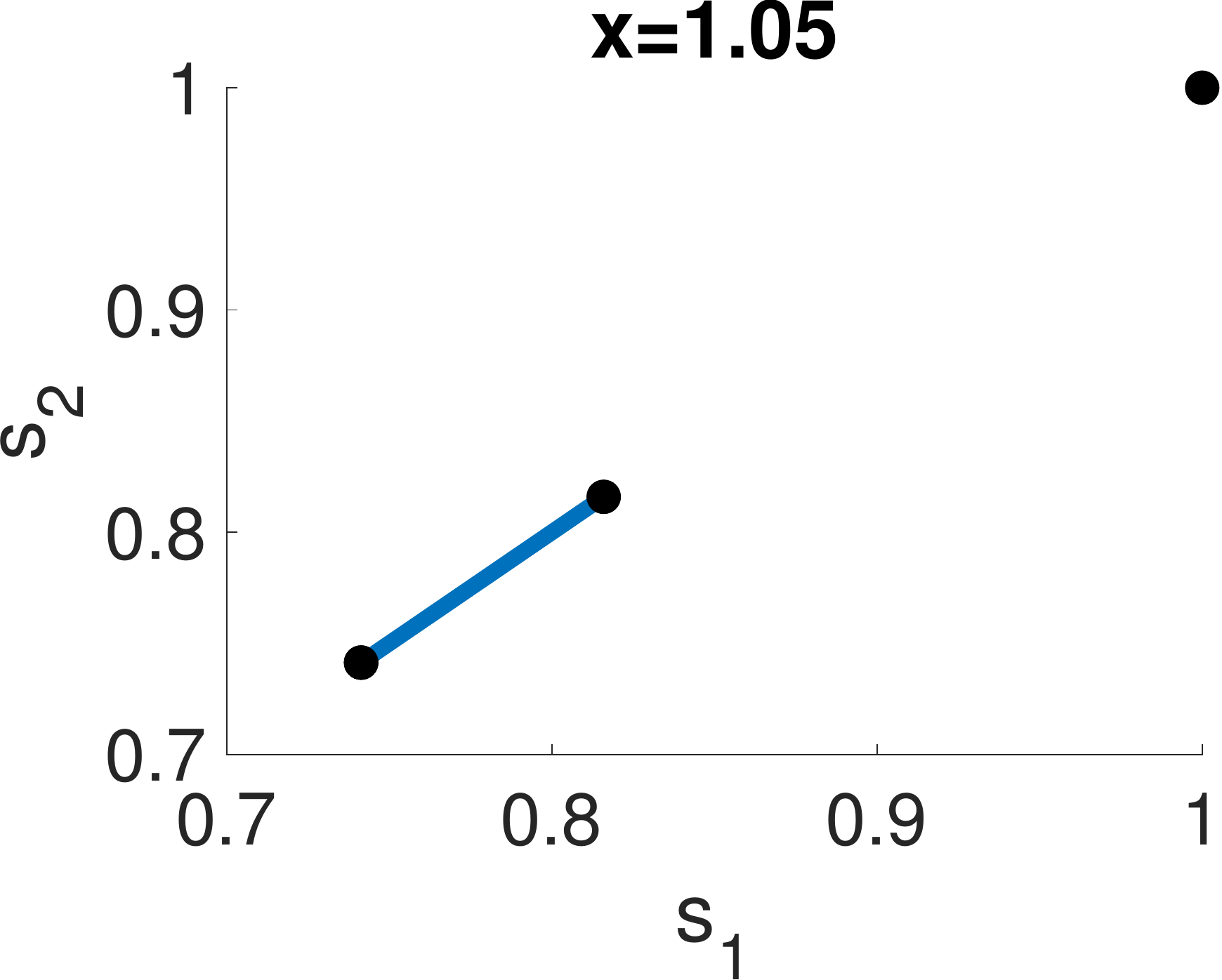}
\includegraphics[width=3.8cm]{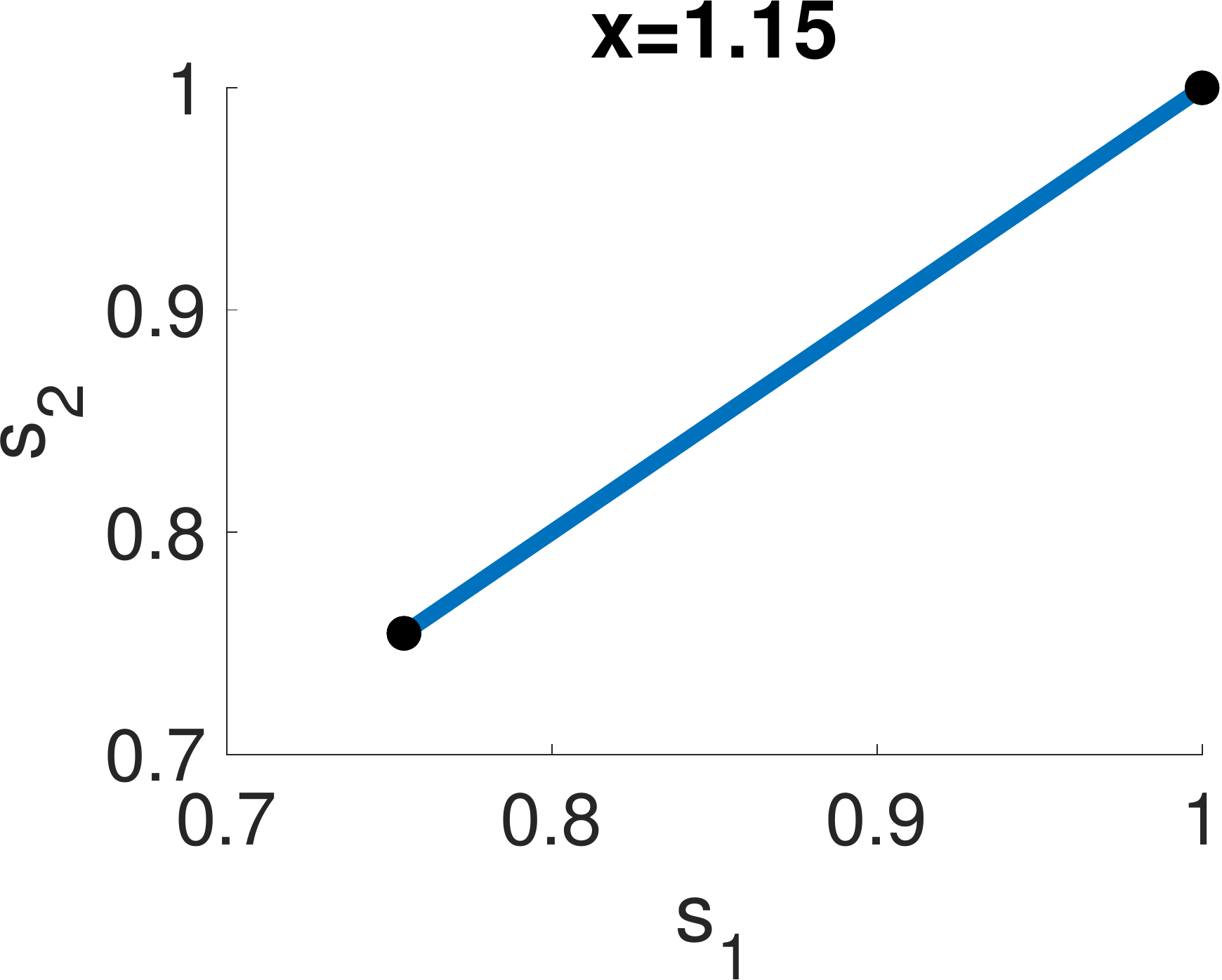}
\vspace{2mm}

\includegraphics[width=3.8cm]{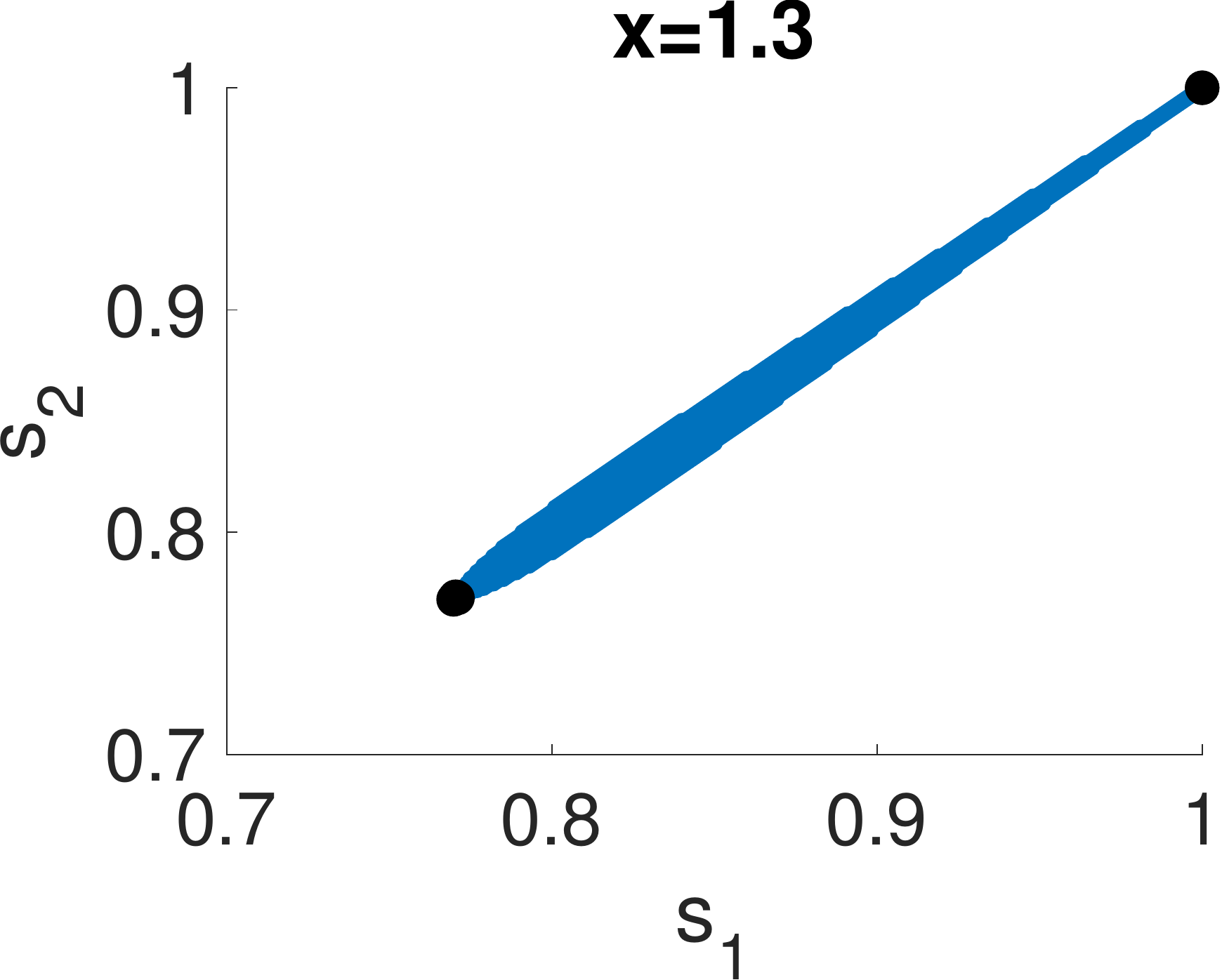}
\includegraphics[width=3.8cm]{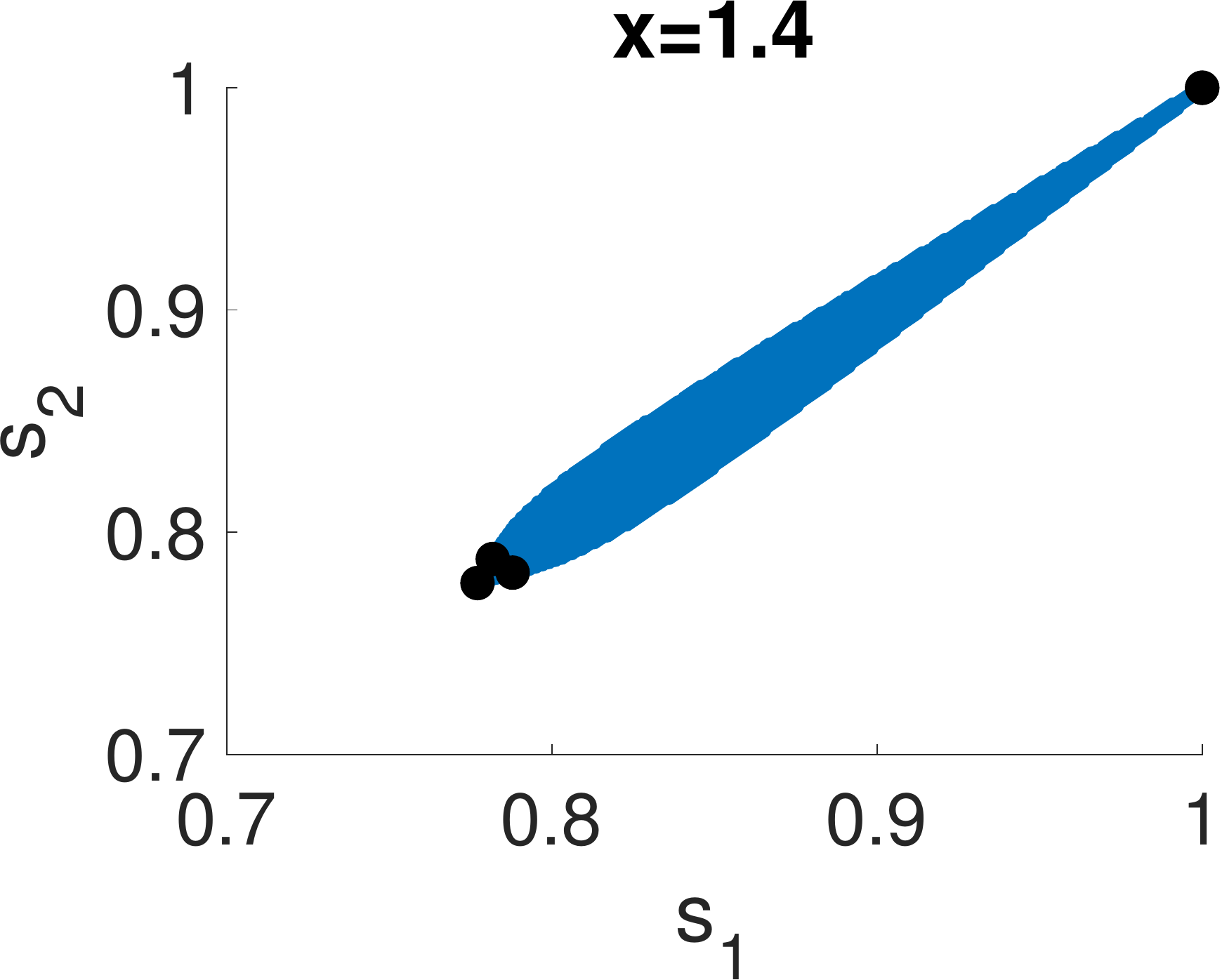}
\includegraphics[width=3.8cm]{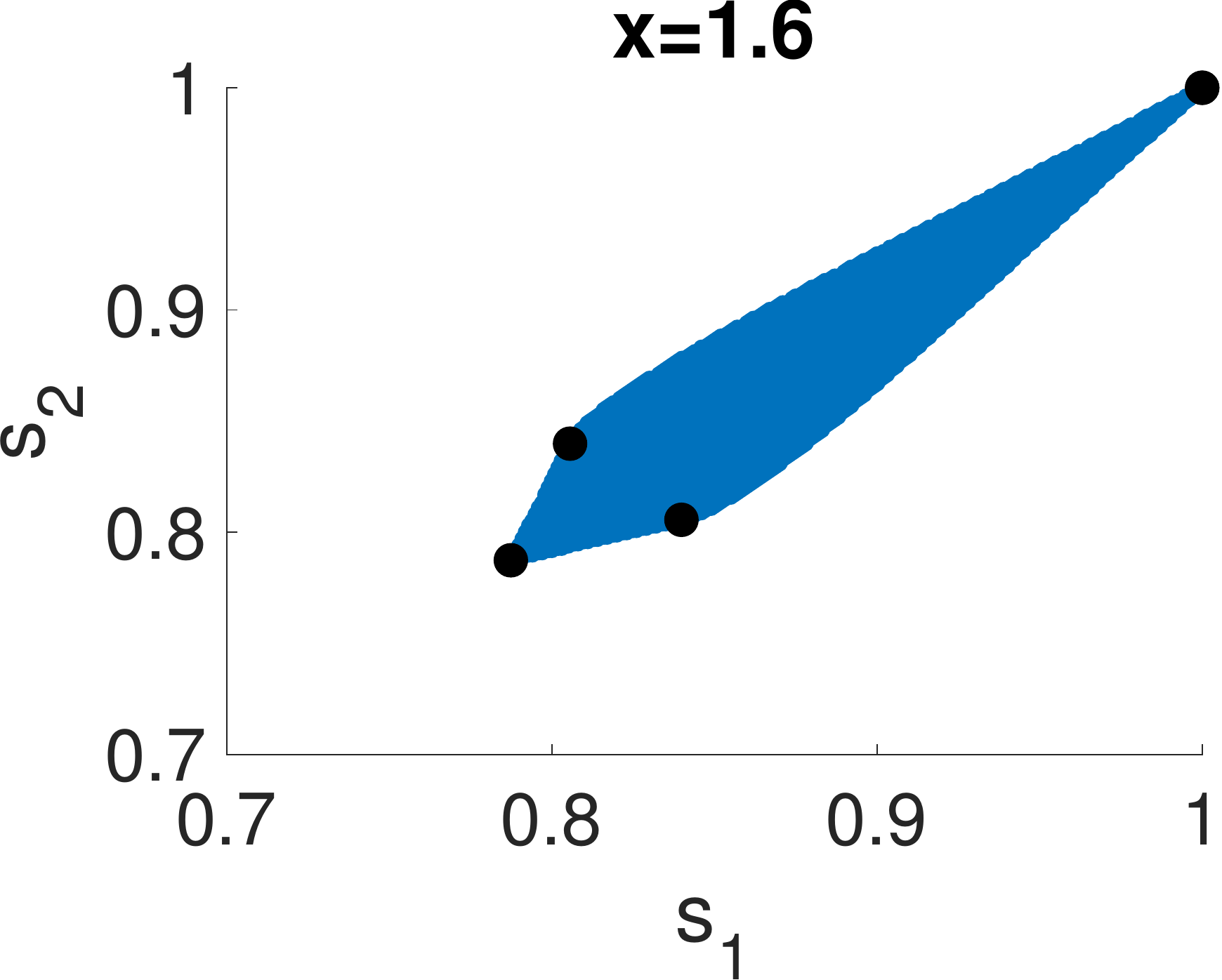}
\vspace{2mm}

\includegraphics[width=3.8cm]{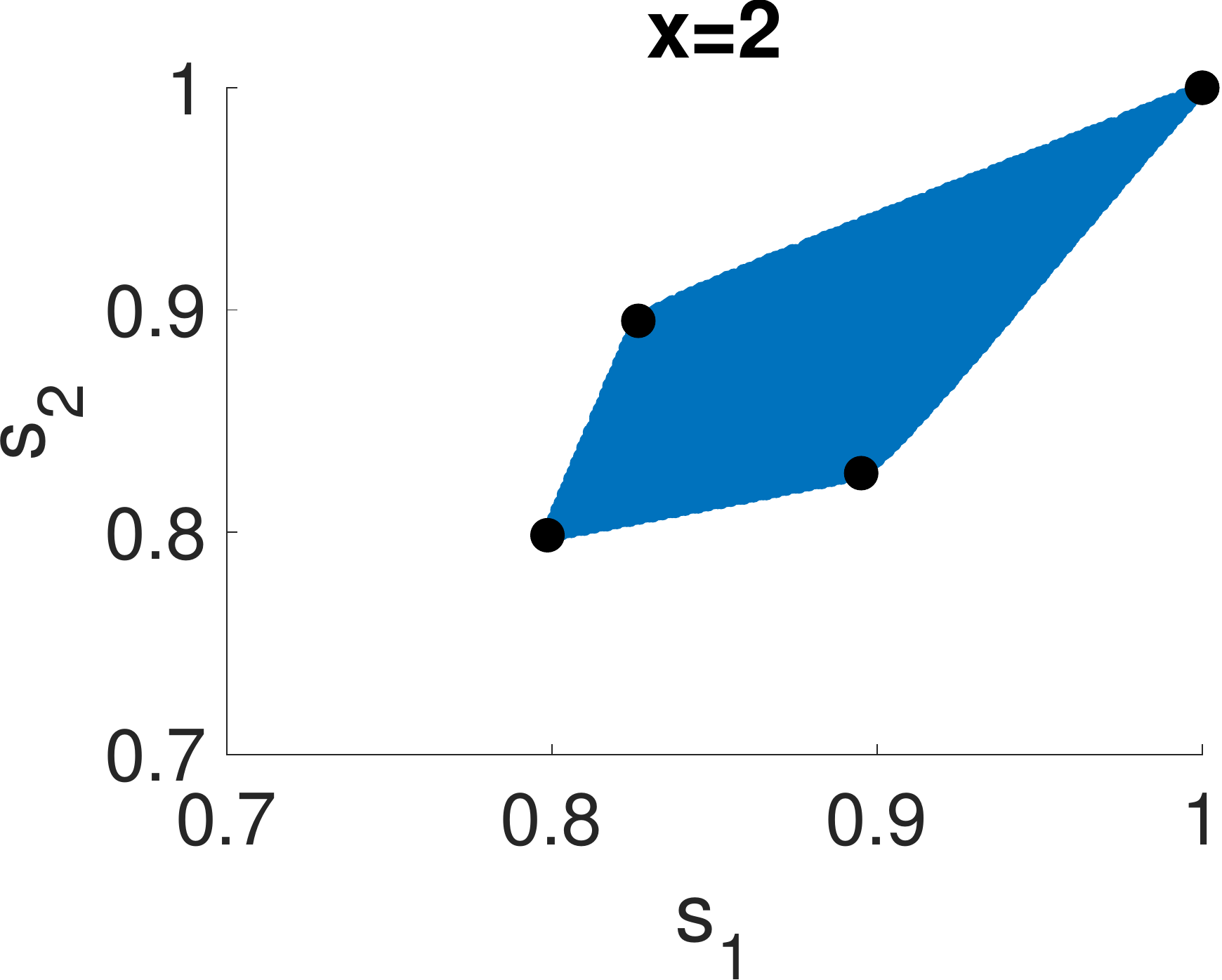}
\includegraphics[width=3.8cm]{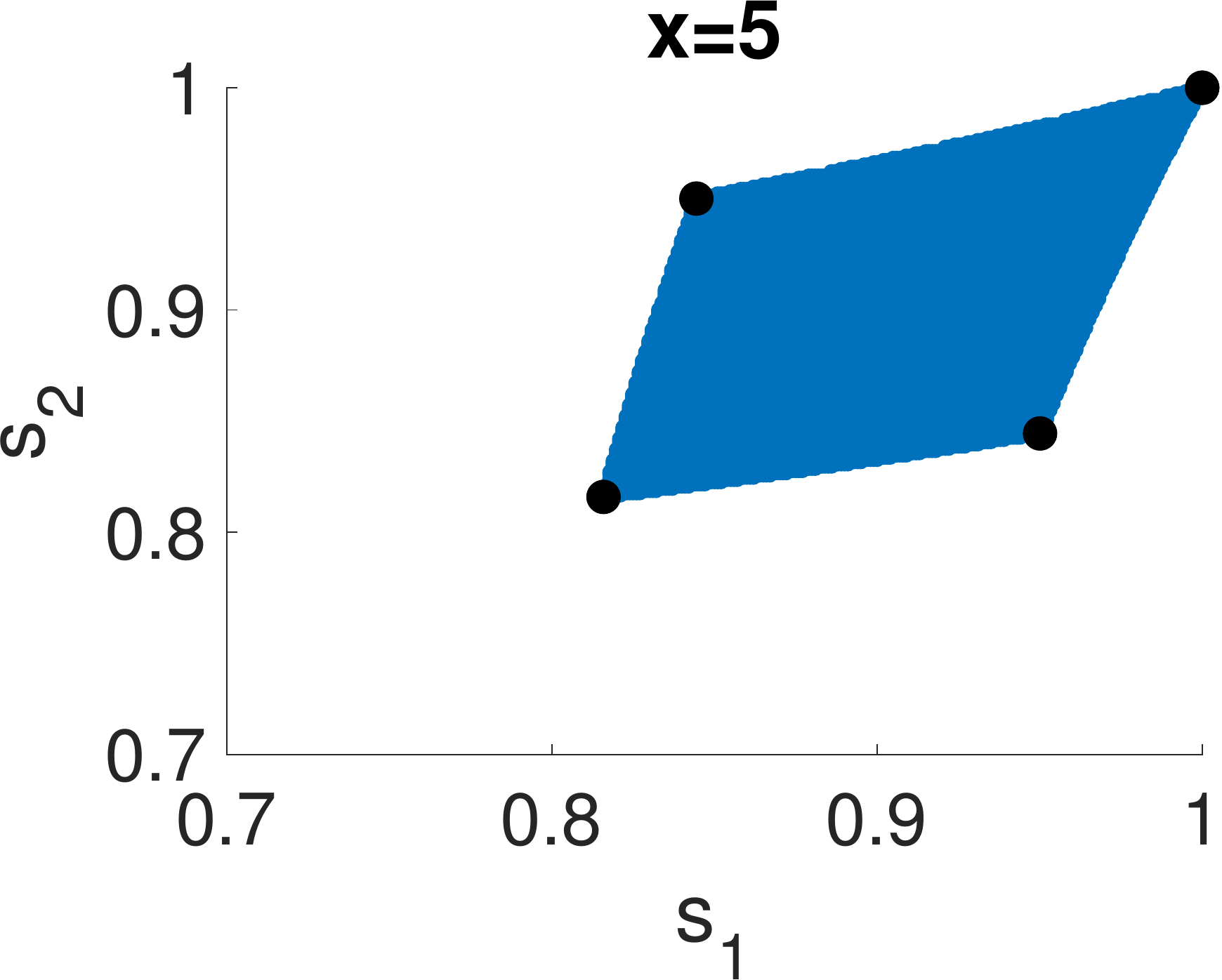}
\includegraphics[width=3.8cm]{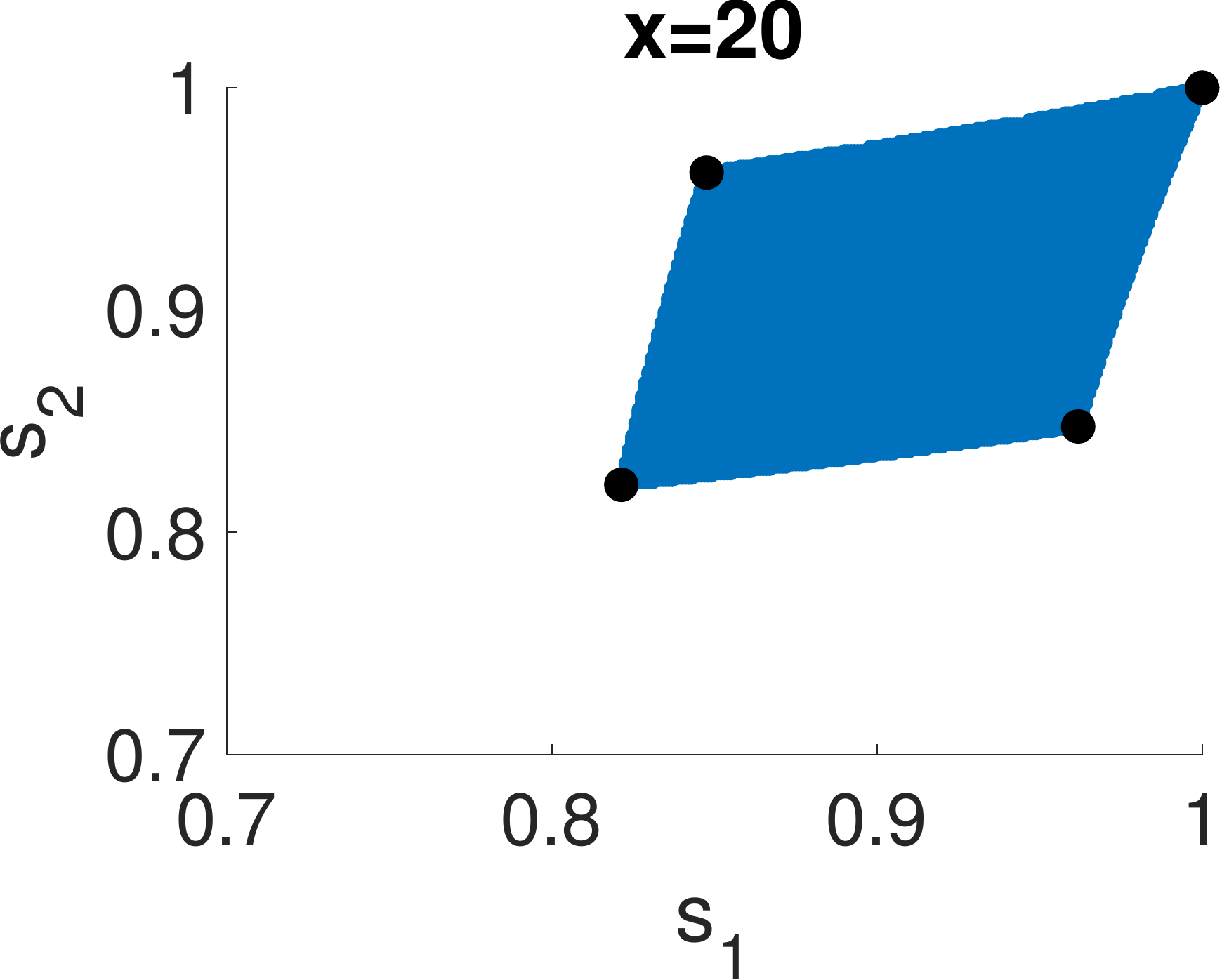}
\caption{ \label{Exam1}Example 1. Top: the extinction probabilities $q_{\langle 0,1 \rangle}$, $q_{\langle 0,1 \rangle}(A_1)$, $q_{\langle 0,1 \rangle}(A_2)$ and $\tilde q_{\langle 0,1 \rangle}$  for $1\leq x\leq 3$. Bottom: the projection set $S_{0}$ for nine values of $x$   (with the shorthand notation $s_i$ for $s_{\langle 0,i\rangle}$, $i=1,2$). }
\end{figure}

\medskip

 We let $a=1/5$, $b=0$, $c=1$, and $y=1/5$, and we study the extinction probabilities for different values of the parameter $x$. In this case, $b+y+2\sqrt{ac}\approx 1.09$ and $\mu \approx 1.38$.  Here $\sup_{i \in A_j } \tilde{q}^{(A_j)}_i<1$ for $j=1,2$, therefore we can use Theorem \ref{COD2} to compute  $\vc q(A_1)$ and $\vc q(A_2)$.
The top graph in Figure~\ref{Exam1} depicts the extinction probabilities $q_{\langle 0,1 \rangle}$, $q_{\langle 0,1 \rangle}(A_1)$, $q_{\langle 0,1 \rangle}(A_2)$ and $\tilde q_{\langle 0,1 \rangle}$  for $1\leq x\leq 3$.  
As Proposition \ref{SDExam} suggests, we observe two phase transitions points, the first at $x=1$, where the number of distinct extinction probability vectors increases from one to two (even if it only becomes clear slightly after $x=1$), and the second at $x=\mu$, where the number of distinct extinction probability vectors increases from two to four. By visual inspection, there exists an  infimum value of $x$ for which $\bm{\tilde{q}} =\bm{1}$. Using Theorem \ref{globpart_crit} we numerically estimate this value to be approximately $1.09$.
 
   The bottom nine graphs in Figure \ref{Exam1} illustrate the set $S_{0}$ ($S$ projected onto level $0$) for nine values of $x$ ranging from $x=1$ to $x=20$. 
The  projected extinction probabilities $(q_{\langle 0,1 \rangle}, q_{\langle 0,2 \rangle })$, $(q_{\langle 0,1 \rangle}(A_1), q_{\langle 0,2 \rangle}(A_1))$, $(q_{\langle 0,1 \rangle}(A_2),$\\ $q_{\langle 0,2 \rangle}(A_2))$, and $(\tilde q_{\langle 0,1 \rangle}, \tilde q_{\langle 0,2 \rangle})$ are marked  by bold discs.
We observe that for small values of $x$ (i.e $x=1,1.05,1.1$) the elements in $S_0$ cling tightly to the straight line of fixed points connecting $\bm{q}$ and $\bm{\tilde q}$ that we identified in Proposition \ref{IsoFP}.
As $x$ increases, the set $S_{0}$ inflates until it visibly contains area. Noticeably, this occurs when $x\leq \mu$ (i.e. $x=1.2$) as well as when $x>\mu$.
  For large values of $x$ (i.e $x=5,20$), the extinction probabilities $(q_{\langle 0,1 \rangle}(A_1), q_{\langle 0,2 \rangle}(A_1))$ and $(q_{\langle 0,1 \rangle}(A_2), q_{\langle 0,2 \rangle}(A_2))$ appear to correspond to vertices in $S_{0}$, while this is less clear for smaller values of $x>\mu$  (i.e. $x=1.4,1.6$).
  
  Due to the  symmetry of the progeny distributions between phases 1 and 2, the level projection sets $S_k$ are symmetric with respect to the diagonal. The next example considers an asymmetric modification of Example 1.

\paragraph{Example 2}

\begin{figure}
\centering
\begin{tikzpicture}[scale=0.85]

\tikzset{vertex/.style = {shape=circle,draw,minimum size=1.3em}}
\tikzset{edge/.style = {->,> = latex'}}

%

\node[vertex, minimum size=.9cm] (0) at  (0,0) {\footnotesize$0 ,1$};
\node[vertex, minimum size=.9cm] (2) at  (3.5,0) {\footnotesize $1, 1 $};
\node[vertex, minimum size=.9cm] (4) at  (7,0) {\footnotesize $2,1$};
\node[vertex, minimum size=.9cm] (6) at  (10.5,0) {\footnotesize $3,1$};

\node[vertex, minimum size=.9cm] (1) at  (0,2.5) {\footnotesize $ 0, 2 $};
\node[vertex, minimum size=.9cm] (3) at  (3.5,2.5) {\footnotesize $1,2$};
\node[vertex, minimum size=.9cm] (5) at  (7,2.5) {\footnotesize $2,2$};
\node[vertex, minimum size=.9cm] (7) at  (10.5,2.5) {\footnotesize $3,2$};

\node at (13,2.5) {\small $\dots$};
\node at (13,0) {\small $\dots$};



%

\draw[edge,above] (0) to[bend left=15] node {\footnotesize $c$ } (2);
\draw[edge,above] (2) to[bend left=15] node {\footnotesize $c $ } (4);
\draw[edge,above] (4) to[bend left=15] node {\footnotesize $c $ } (6);

\draw[edge,below] (2) to[bend left=15] node {\footnotesize $a$ } (0);
\draw[edge,below] (4) to[bend left=15] node {\footnotesize $a$ } (2);
\draw[edge,below] (6) to[bend left=15] node {\footnotesize $a$ } (4);

\draw[edge,above] (1) to[bend left=15] node {\footnotesize $c$ } (3);
\draw[edge,above] (3) to[bend left=15] node {\footnotesize $c$} (5);
\draw[edge,above] (5) to[bend left=15] node {\footnotesize $c$} (7);

\draw[edge,below] (3) to[bend left=15] node {\footnotesize $a$ } (1);
\draw[edge,below] (5) to[bend left=15] node {\footnotesize $a$} (3);
\draw[edge,below] (7) to[bend left=15] node {\footnotesize $a$} (5);


\draw[edge,left] (0) to[bend left=20] node {\footnotesize $b$ } (1);

\draw[edge,right] (1) to[bend left=20] node {\footnotesize $b$ } (0);
\draw[edge,right] (3) to[bend left=0] node {\footnotesize $b$} (2);
\draw[edge,right] (5) to[bend left=0] node {\footnotesize $b$} (4);
\draw[edge,right] (7) to[bend left=0] node {\footnotesize $b$} (6);

\draw[edge,below] (2) to [out=-110,in=-70,looseness=4.5] node {\footnotesize $b$} (2);
\draw[edge,below] (4) to [out=-110,in=-70,looseness=4.5] node {\footnotesize $b$} (4);
\draw[edge,below] (6) to [out=-110,in=-70,looseness=4.5] node {\footnotesize $b$} (6);


\draw[above] (7) to[bend left=7.5, pos=1] node {\footnotesize $c$} (12.25,2.825);
\draw[edge,below] (12.25,2.175) to[bend left=8,pos=0] node {\footnotesize $a$}  (7);

\draw[above] (6) to[bend left=7.5, pos=1] node {\footnotesize $c$} (12.25,0.325);
\draw[edge,below] (12.25,-0.325) to[bend left=8,pos=0] node {\footnotesize $a$}  (6);

%

\end{tikzpicture}
\caption{\label{Exam2MPRG}Mean progeny representation graph corresponding to Example 2.}
\end{figure}
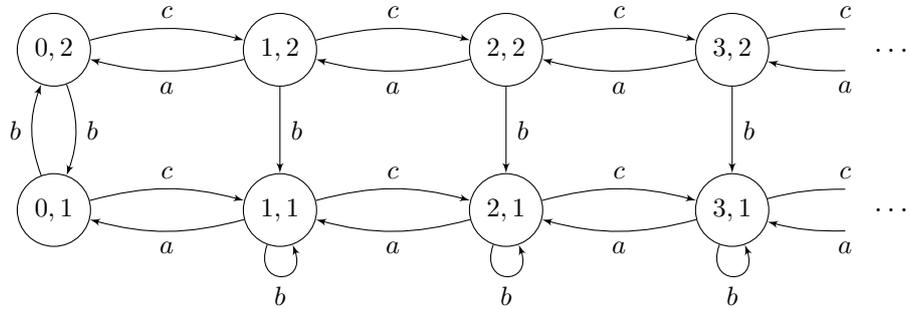

\begin{figure}
\centering
\includegraphics[width=10cm, angle = -90]{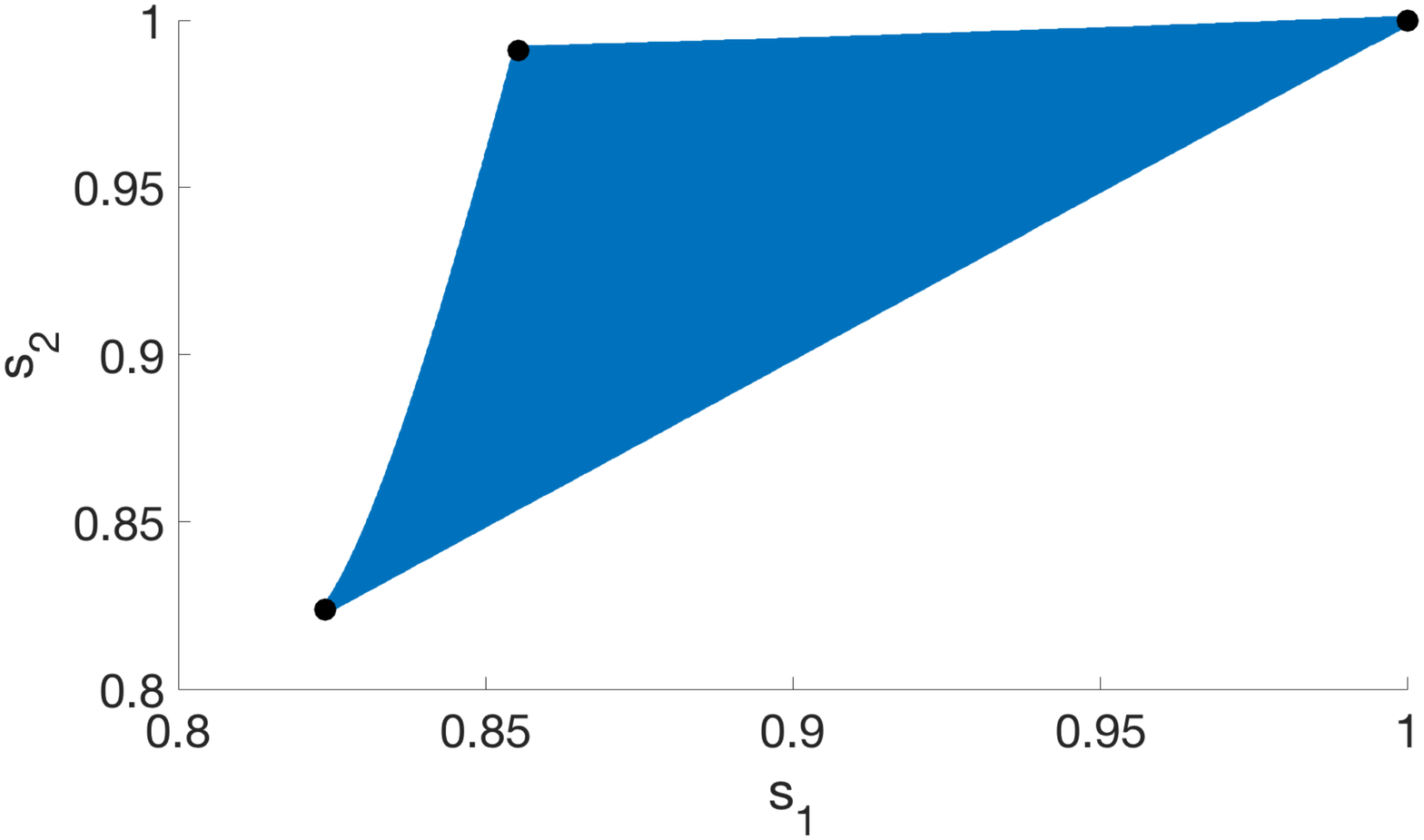}
\caption{ \label{Exam2fp}The projection set $S_{0}$ corresponding to Example 2. }
\end{figure}

We modify Example 1 such that the mean progeny representation graph becomes as shown in Figure \ref{Exam2MPRG}. In Figure \ref{Exam2fp} we plot the set $S_0$ for $a=1/5$, $b=1/20$, and $c=1$, and observe that there are only three distinct extinction probability vectors, $\bm{q}$, $\bm{q}(A_2)$ and $\bm{ \tilde q}$. Indeed, by   Corollary~\ref{InfiniteGlobal}, $\bm{q}(A_1) = \bm{q}$. 
Despite the lack of symmetry, this branching process is still locally isomorphic to a single phase LHBP. Thus, by Proposition~\ref{IsoFP}, the set $S_0$ still contains the linear segment that connects the global and partial extinction probabilities, $(q_{\langle 0,1 \rangle} , q_{\langle 0 ,2 \rangle})$ and $(\tilde q_{\langle 0,1 \rangle} , \tilde q_{\langle 0 ,2 \rangle})$. 
On inspection of Figure \ref{Exam2fp} we see that this linear segment now sits on the boundary of $S_0$. 

\medskip

Examples 1 and 2 motivate several questions, which to our knowledge remain open. In particular: 
\begin{itemize}
\item[\emph{(i)}] We have only focused on sets $A \in \sigma(A_1, \dots, A_d)$, leading to a maximum of $2^d$ potentially distinct extinction probability vectors $\vc q(A)$; we may then ask if more than $2^d$ distinct extinction probability vectors can exist in an irreducible block LHBP when we consider \emph{any} $A \subseteq \mathcal{X}_d$. 
\item[\emph{(ii)}] Given the set $S$, we may question whether it is possible to identify which elements correspond to extinction probability vectors.
\end{itemize}

We now propose an answer to \emph{(ii)}, which we suggest applies to any irreducible multitype Galton-Watson branching processes with countably many types:

\begin{conjecture}\label{ConjectPH}
If $\bm{q}=\bm{\tilde q}$ then $S=\{ \bm{q}, \bm{1} \}$, whereas if $\bm{q}<\bm{\tilde q}$ then $S$ contains a continuum of elements, whose minimum is $\bm{q}$,  and whose maximum is $\bm{\tilde q}$. In addition, the boundary of any projection set is differentiable everywhere except at points that correspond to an extinction probability vector $\bm{q}(A)$ for some \mbox{$A \subseteq \mathcal{X}_d$}.
\end{conjecture}

\section*{Acknowledgements}
The authors would like to acknowledge the support of the Australian Research Council (ARC) through the Centre of Excellence for the Mathematical and Statistical Frontiers (ACEMS). Sophie Hautphenne would further like to thank the ARC for support through Discovery Early Career Researcher Award DE150101044.

\end{document}